%% file: main.tex
\makeatletter
\let\SIAMkernelrefstepcounter\refstepcounter
\AddToHook{package/hyperref/before}{\let\refstepcounter\SIAMkernelrefstepcounter}
\makeatother
\documentclass[final,onefignum,onetabnum]{siamart220329}
\usepackage[T1]{fontenc}
\usepackage[utf8]{inputenc}
\usepackage{amsmath,amssymb,amsfonts,booktabs,graphicx}
\usepackage{microtype}
\usepackage{placeins}
\usepackage{cleveref}
\usepackage{algorithmicx}
\usepackage{algpseudocode}
\usepackage{tikz}
\usetikzlibrary{arrows.meta,shapes.geometric}
\newsiamremark{remark}{Remark}
\newsiamremark{example}{Example}
\newcommand{\HH}{\mathcal H}

\newcommand{\Id}{\operatorname{Id}}
\newcommand{\Fix}{\operatorname{Fix}}
\newcommand{\zer}{\operatorname{zer}}

\newcommand{\ip}[2]{\langle #1,#2\rangle}
\newcommand{\norm}[1]{\lVert #1\rVert}
\headers{Second-order Accelerated Halpern-type Method}{Shuchang Zhang, Xianchen Zhou, and Jianping Zhang}
\title{SAHM: Second-order Accelerated Halpern-type Method with $\lowercase{o}(1/\lowercase{k})$ Fixed-point Residual Rate}
\author{Shuchang Zhang\thanks{School of Mathematics and Computational Science,
Xiangtan University, Xiangtan, China.}
\and Xianchen Zhou\thanks{College of Science, National University of Defense Technology, China.}
\and Jianping Zhang\thanks{School of Mathematics and Computational Science,
Xiangtan University, Xiangtan, China. Corresponding author.}}
\hypersetup{pdftitle={SAHM: Second-order Accelerated Halpern-type Method with o(1/k) Residual Decay and Convergence Bias Quantification},
pdfauthor={Shuchang Zhang and Jianping Zhang}}
\begin{document}
\raggedbottom
\widowpenalty=10000
\clubpenalty=10000
\maketitle
\begin{abstract}
Accelerating fixed-point residual of Halpern iteration is a hot topic, and second-order dynamics can generate accelerated fast Krasnosel'ski\u{\i}--Mann (KM)  iteration with $\lowercase{o}(1/k)$ rate of fixed-point residual. However, the classic Halpern iteration generally admits the $ O(1/k)$ fixed-point residual. The fast KM iteration converges weakly in the Hilbert space, and the limit point lacks a quantitative characterization. Building on Tikhonov regularization and second-order dynamics, we propose the Second-order Accelerated Halpern-type Method (SAHM), which not only converges strongly, but also achieves the $\lowercase{o}(1/k)$ rate of the fixed-point residual. We develop a quantitative analysis of convergence bias, projection geometry and variational defects quantify the discrepancy between their limit and the projection of the initial point. The weak limit points of KM iteration and its accelerated variant may distinct from the projections of the initial point. Numerical Experiments show the $\lowercase{o}(1/k)$ rate and the competitiveness of SAHM.
\end{abstract}
\begin{keywords}
Halpern iteration, second-order dynamics, strong convergence, fixed-point residual, Tikhonov regularization, Fast Krasnosel'ski\u{\i}--Mann method
\end{keywords}
\begin{MSCcodes}
47H09, 47H10, 47J25, 65J15, 65K15
\end{MSCcodes}
\input{sections/introduction}
\input{sections/preliminaries}
\input{sections/continuous}
\input{sections/discrete}
\input{sections/selection}
\input{sections/experiments}

\FloatBarrier
\input{sections/discussion}
\appendix
\input{sections/energy_appendix}
\bibliographystyle{siamplain}
\bibliography{references}
\end{document}

%% file: sections/introduction.tex
\section{Introduction}\label{sec:intro}
Let $T$ be a nonexpansive mapping on a real Hilbert space $\HH$ with $Z=\Fix T\ne\varnothing$, and let $p=P_Zu$ for a prescribed anchor $u$. The Krasnosel'ski\u{\i}--Mann (KM) iteration
\begin{equation}\label{intro:km}
 x_{k+1}=(1-\lambda_k)x_k+\lambda_k Tx_k
\end{equation}
originates in \cite{Krasnoselskii,Mann}. Under standard relaxation
conditions it converges weakly to a fixed point. 
KM is Fej\'er monotone with respect to $Z$, which provides a precise sense in which it stays near its initial point. It nevertheless need not converge to $P_Zx_0$, even for compositions of simple projections.

Halpern's iteration \cite{Halpern} introduces a fixed anchor:
\begin{equation}\label{intro:halpern}
 x_{k+1}=\lambda_k u+(1-\lambda_k)Tx_k.
\end{equation}
Under suitable conditions on $\{\lambda_k\}$ it converges strongly to $P_Zu$~\cite{Wittmann,Xu2004}. It has the strong convergence advantage than KM iteration in Hilbert space~\cite{He-Xu-Dong-Mei}. This observation shows that the Halpern iteration is an implicit regularized method, and it has a distinct convergence bias compared with the KM iteration. Viscosity approximation, introduced by Moudafi \cite{Moudafi2000}
and developed by Xu \cite{Xu2004}, replaces the fixed anchor by a
contraction and selects the solution of a variational inequality on
$Z$. Sabach and Shtern \cite{Sabach-Shtern} obtained an $O(k^{-1})$
residual estimate through the sequential averaging method.
Lieder \cite{Lieder} obtained the tight estimate
$\norm{x_k-Tx_k}\le2\norm{u-p}/(k+1)$ for $x_0=u$ and
$\lambda_k=1/(k+2)$. This is an open-loop rule: the weight $\lambda_k$ depends on the iteration number, rather than history observations. He et al. developed adaptive anchoring rules with strong convergence, refined residual estimates, and improved numerical performance~\cite{He-Xu-Dong-Mei}.

Anchoring also provides a design principle for monotone inclusions
and minimax problems. Yoon and Ryu \cite{Yoon-Ryu:21} combined it with
extragradient steps to obtain $O(k^{-2})$ decay of the squared
gradient norm. Related accelerated splitting and extragradient
constructions were developed by Tran-Dinh and Luo
\cite{Tran-Dinh-Luo} and Lee and Kim \cite{Lee-Kim}. The merging-path analysis of Yoon
and Ryu \cite{Yoon-Ryu} clarifies why several accelerated trajectories
approach one another. Diakonikolas \cite{Diakonikolas2020} developed
near-optimal and parameter-free methods for monotone inclusion and
variational inequalities. Actually, the discrete KM, Halpern, and extragradient methods can be derived from a continuous-time dynamical system. Furthermore,Tikhonov regularization and second-order dynamics provide a novel perspective for the design and acceleration of Halpern-type methods. 
\subsection{Tikhonov regularization and second-order dynamics}
Tikhonov regularization is an established route to strong convergence with a particular solution \cite{Bot-Csetnek-Meier,Bot-Meier,Bot-Chenchene-Flexible,Bot-Nguyen-Tikhonov}. Monotone operator inclusion problems can be related to dynamical systems. The inclusion problem $A(x) =0$ can be considered from the following dynamical system:
\begin{equation}\label{eq:firsrt-order-dynamics}
 \dot x(t)+A(x(t))=0, 
\end{equation}
where $t\geq t_0 (t_0>0)$. However, it is a well-known fact that the trajectories of~\eqref{eq:firsrt-order-dynamics} can only be expected to converge weakly and ergodically to a zero of $A$~\cite{baillon1976remark}. It considers the following regularized system:
\begin{equation}\label{eq:firsrt-order-regularized dynamics}
 \dot x(t)+A(x(t))+\varepsilon(t)(x(t)-u)=0,
\end{equation}
where $u\in \mathcal{H}$ is the anchor term. We can replace $A$ with  $V = A+\varepsilon(I-u)$ in \eqref{eq:firsrt-order-dynamics}, which yields \eqref{eq:firsrt-order-regularized dynamics}. If $\varepsilon$ satisfies certain special properties, the trajectory generated by \eqref{eq:firsrt-order-regularized dynamics} may exhibit strong convergence, providing a continuous dynamical-system perspective for studying the discrete Halpern iteration. Bo\c{t} and Nguyen \cite{Bot-Nguyen-Tikhonov} analyzed first-order regularized monotone flows and their links to Halpern iterations and dynamics systems. Bo\c{t} and Chenchene \cite{Bot-Chenchene-Flexible} obtained strong selection and $O(k^{-1})$ residual decay for an extragradient method with flexible
anchoring. Bo\c{t}, Chenchene, and Fadili
\cite{Bot-Chenchene-Fadili} further developed the relation between inertial and anchoring mechanisms, with overrelaxation and preconditioning. Accelerated diagonal methods for bilevel optimization \cite{Bot-Chenchene-Csetnek-Hulett} provide another related setting in which regularization decay interacts with geometric assumptions.

Second-order dynamics have been used to accelerate the fixed-point residual.
Bo\c{t}, Csetnek, and Nguyen \cite{Bot-Csetnek-Nguyen} developed the
Fast OGDA dynamics and discrete methods for monotone equations. For nonexpansive mappings, Bo\c{t} and Nguyen \cite{Bot-Nguyen-FastKM} obtained a Fast KM iteration with weak convergence and fixed-point residual $o(k^{-1})$. A central continuous model, with $A=I-T$, is
\begin{equation}\label{intro:fastflow}
 \ddot x+\frac{\alpha}{t}\dot x+\frac{\mathrm{d}}{\mathrm{d}t}A(x)
                         +\frac{\alpha}{2t}A(x)=0,\qquad \alpha>2.
\end{equation}
The extra residual-difference term in its discretization is essential.
It cannot be recovered by changing only the relaxation weight in \eqref{intro:km}. 

However, existing Halpern iterations generally achieve a $\mathcal{O}(k^{-1})$ fixed-point residual rate, and obtaining $o(k^{-1})$ rate remains open. Although the Fast KM iteration can attain an $o(k^{-1})$ rate, its iterates converge only weakly, and there is no quantitative characterization of the weak-limit points of the KM sequence: their relation to the projection of the initial point is unclear, and it is unknown whether Fast KM exhibits the solution-selection property of Halpern iterations. Two objectives arise when computing a fixed point: reducing the residual $\norm{x-Tx}$ and selecting a particular
element of $Z$. These objectives are different whenever $Z$ contains more than one point. A small residual certifies approximate feasibility; it does not identify the feasible point nearest to a given initial point.  
\subsection{Contributions}
Motivated by Tikhonov regularization and second-order dynamics, we regularize the operator in \eqref{intro:fastflow} by replacing $A$ with $V = A+\varepsilon(I-u)$, then we can derive the explicit second-order dynamics with Halpern-type anchoring:
\begin{equation}
 \ddot x+\left(\frac{\alpha}{t}+\varepsilon\right)\dot x
 +\frac{\mathrm{d}}{\mathrm{d}t}A(x)+\frac{\alpha}{2t}A(x)
 +\left(\dot\varepsilon+\frac{\alpha}{2t}\varepsilon\right)(x-u)=0.
 \label{ct:expanded-flow}
\end{equation}
Existing works on the development of Halpern-type methods have focused on the rate of the fixed-point residual $\left\| Tx^k-x^k \right\| $. The classical Halpern iteration generally produces the ${O}(k^{-1})$ original residual on this fixed problem. Our aim is the simultaneous strong convergence property and pointwise $o(k^{-1})$ residual decay by combining a fixed anchor with inertia and a residual-difference correction. The proposed \emph{Second-order accelerated Halpern-type method (SAHM)} is summarized in \cref{alg:sahm}. By reusing the previous operator value, each iteration requires only one new evaluation of $A$.
\begin{algorithm}[t]
\caption{Second-order accelerated Halpern-type method (SAHM)}
\label{alg:sahm}
\begin{algorithmic}[1]
\Require Initial points $x_0,x_1$, anchor $u$, parameters
$\alpha,s$, and anchoring schedule $(\varepsilon_k)_{k\ge0}$.
\State Set $g_0=A(x_0)$.
\For{$k=1,2,\ldots$}
    \State Evaluate $g_k=A(x_k)$.
    \State Set
    \[
    \begin{aligned}
    &\beta_k=\frac{k(1-s\varepsilon_{k-1})}{k+\alpha}, \tau_k=\frac{s\alpha}{2(k+\alpha)},\\
    &\sigma_k=\frac{sk}{k+\alpha},\omega_k=\frac{s}{k+\alpha}\left[\frac{\alpha}{2}\varepsilon_k
    +k(\varepsilon_k-\varepsilon_{k-1})\right].
    \end{aligned}
    \]
    \State Update
    \begin{equation*}
    x_{k+1}= x_k+\beta_k(x_k-x_{k-1})-\tau_k g_k-\sigma_k(g_k-g_{k-1})-\omega_k(x_k-u).
    \end{equation*}
\EndFor
\end{algorithmic}
\end{algorithm}
Here, Halpern-type methods refer to methods converging strongly to $P_Zu$, while acceleration refers to the fixed-point residual $\left\| Tx^k-x^k \right\| $ decay. The method retains the inertial and residual-difference structure of Fast KM \cite{Bot-Nguyen-FastKM}; Tikhonov regularization supplies the anchoring mechanism. 

The main contributions of this paper are as follows:
\begin{itemize}
\item \Cref{ct:main} establishes strong convergence to $P_Zu$ and
$o(t^{-1})$ decay of the original residual and velocity for a second-order
flow with a monotone Lipschitz operator and an independently prescribed
anchor. A cancellation-based correction to the Fast OGDA energy
motivates the discrete construction.

\item \Cref{dis:main} establishes strong convergence to $P_Zu$ and $o(k^{-1})$ decay of the original residual and increments for \cref{alg:sahm}. Its independent Fast KM energy argument controls the operator-history error through anchoring, allowing $\varepsilon_k=c/[(k+\kappa)\log^q(k+\kappa)]$ for all $0<q\le1$, even when the resulting energy errors are nonsummable.

\item We quantify the discrepancy between unregularized limits and $P_Zu$ through KM dissipation, an exact Fast KM path identity, affine-set formulas, and a nonaffine example with certified distinct limits. Projection geometry and a normal-cone defect bound the secondary objective and selection error. KM and Fast KM approach fixed points distinct from $P_Zx_0$, which have the persistent selection bias.

\item Experiments on polyhedral and curved sets assess the fixed-point residual, the distance between the iterates and the true solution, and discrepancy between the limit point and the projection of the initial point. We also compare EG~\cite{Korpelevich1976}, EAG-C and EAG-V~\cite{Yoon-Ryu:21}, FEG~\cite{Lee-Kim}, and
the anchored Popov method (APV)~\cite{Tran-Dinh-Luo}, and the logarithmic SAHM (Log-SAHM) on the same problems with identical initial points and anchors. The proposed SAHM achieves state-of-the-art oracle efficiency for attaining joint residual and projection accuracy on the tested problems.
\end{itemize}

After the preliminary comparison principles in \cref{sec:prelim},
\cref{ct:section} develops the continuous dynamics and its convergence theory. \Cref{sec:discrete} then gives the discrete SAHM algorithm and its convergence analysis, followed by the solution selection analysis in \cref{sel:section} and the experiments in \cref{num:section}. The appendix contains the explicit energy equivalence and forward estimates used in the discrete proof.

%% file: sections/preliminaries.tex
\section{Preliminaries}\label{sec:prelim}
Throughout, $\HH$ is a real Hilbert space. A mapping $A:\HH\to\HH$ is
$\rho$-cocoercive, with $\rho>0$, if
\[
 \ip{Ax-Ay}{x-y}\ge \rho\norm{Ax-Ay}^2 \quad(x,y\in\HH).
\]
Such an operator is monotone and $1/\rho$-Lipschitz. If $T$ is
nonexpansive, then $A=I-T$ is $1/2$-cocoercive; if $T$ is
$\theta$-averaged, the constant improves to $1/(2\theta)$. For a continuous monotone operator defined on all of $\HH$,
the zero set $Z=\zer A$, when nonempty, is closed and convex.
Its graph is closed for weak convergence of the argument and strong
convergence of the value. These standard facts, and the corresponding
demiclosedness principle for $I-T$, can be found in
\cite{Bauschke-Combettes:book}.
We fix an anchor $u\in\HH$ and write $p=P_Zu$. Then
\begin{equation}\label{pre:projection}
 \ip{z-p}{u-p}\le0,\qquad\forall z\in Z.
\end{equation}

\begin{lemma}[{\cite{Bot-Nguyen-Tikhonov}}]
\label{pre:continuous-contraction}
Let $t_0\ge0$, $\eta>0$, and let
$\varepsilon:[t_0,\infty)\to[0,\infty)$ be locally integrable with
$\int_{t_0}^{\infty}\varepsilon(t)\,dt=\infty$.
Let $F:[t_0,\infty)\to[0,\infty)$ be locally absolutely continuous and
$b:[t_0,\infty)\to\mathbb R$ be measurable.
Suppose that, for almost every sufficiently large $t$,
\[
 F'(t)\le-\eta\varepsilon(t)F(t)+\varepsilon(t)b(t),
 \qquad \limsup_{t\to\infty}b(t)\le0.
\]
Then $F(t)\to0$ as $t\to\infty$.
\end{lemma}

\begin{lemma}[\cite{Xu2002}]\label{pre:contraction}
Let $E_k\ge0$, $\eta>0$, and $\varepsilon_k>0$ satisfy
$\varepsilon_k\to0$ and $\sum_k\varepsilon_k=\infty$.
If, for all sufficiently large $k$,
\[
 E_{k+1}\le(1-\eta\varepsilon_k)E_k+\varepsilon_k b_k,
 \qquad \limsup_{k\to\infty}b_k\le0,
\]
then $E_k\to0$.
\end{lemma}

%% file: sections/continuous.tex
\section{A second-order flow with Halpern-type anchoring}
\label{ct:section}

We first establish strong selection and an $o(t^{-1})$ residual for a
second-order flow with a prescribed fixed anchor. This analysis applies to monotone
Lipschitz operators. The flow motivates the discrete method developed
in \cref{sec:discrete}, whose convergence proof uses cocoercivity to
control the mismatch between consecutive regularized operators.

Let $\mathcal H$ be a real Hilbert space, let
$A\colon\mathcal H\to\mathcal H$ be monotone and globally Lipschitz, and assume
that $Z:=\operatorname{zer}A$ is nonempty.  The set $Z$ is closed and convex.
Fix an anchor $u\in\mathcal H$ and write $p=P_Zu$.  For $t_0>0$ and a positive
function $\varepsilon\in C^1([t_0,\infty))$, define
\begin{equation*}
 V(t,x)=A(x)+\varepsilon(t)(x-u).
\end{equation*}
We consider the initial-value problem
\begin{equation}
\begin{cases}
& \ddot x(t)+\frac{\alpha}{t}\dot x(t)
 +\frac{\mathrm{d}}{\mathrm{d}t}V(t,x(t))
 +\frac{\alpha}{2t}V(t,x(t))=0, \\
& x(t_0)=x_0, \dot x(t_0)=v_0,
\end{cases}
 \label{ct:flow}
\end{equation}
where $\alpha>2$ and $x_0,v_0\in\mathcal H$ are arbitrary.
The derivative of $V(t,x(t))$ is the total derivative.  In particular,
\eqref{ct:flow} expands as~\eqref{ct:expanded-flow}. The term $\dot\varepsilon(x-u)$ is part of the proposed system.

\begin{proposition}[First-order representation and well-posedness]
\label{ct:well-posedness}
Under the preceding assumptions, \eqref{ct:flow} has a unique global solution
$x\in C^1([t_0,\infty);\mathcal H)$ with locally absolutely continuous
derivative.  The equation holds almost everywhere.  More precisely, setting
$a_0=\alpha/2>1$ and $b_0=a_0-1$, it is equivalent to
\begin{equation}
 \dot x=\frac{a_0}{t}(w-x)-A(x)-\varepsilon(x-u),
 \qquad
 \dot w=\frac{b_0}{t}(x-w),
 \label{ct:first-order}
\end{equation}
with initial conditions
\begin{equation}
 x(t_0)=x_0,\qquad
 w(t_0)=x_0+\frac{t_0}{a_0}
 \bigl(v_0+A(x_0)+\varepsilon(t_0)(x_0-u)\bigr).
 \label{ct:first-order-data}
\end{equation}
\end{proposition}

\begin{proof}
The vector field in \eqref{ct:first-order} is continuous and globally
Lipschitz in $(x,w)$, uniformly on compact time intervals, with linear
growth there.  It therefore has a unique global $C^1$ solution.  Moreover,
\begin{equation}
 w=x+\frac{t}{a_0}\bigl(\dot x+V(t,x)\bigr).
 \label{ct:auxiliary-variable}
\end{equation}
$A(x(\cdot))$ and $\dot x$ are locally absolutely continuous by the
Lipschitz property of $A$ and \eqref{ct:first-order}.  Differentiating
\eqref{ct:auxiliary-variable} almost everywhere proves the equivalence
with \eqref{ct:flow}; its initial value gives \eqref{ct:first-order-data}.
\end{proof}

\begin{theorem}[Strong selection and a vanishing scaled residual]
\label{ct:main}
Suppose that $A\colon\mathcal H\to\mathcal H$ is monotone and globally
Lipschitz, $Z=\operatorname{zer}A\ne\varnothing$, and $\alpha>2$.
Fix $u\in\mathcal H$ and $t_0>0$, and let
$\varepsilon\in C^1([t_0,\infty))$ be positive and satisfy
\begin{equation}
 \int_{t_0}^{\infty}\varepsilon(t)\,\mathrm{d}t=\infty,
 \qquad
 t\varepsilon(t)\longrightarrow0,
 \qquad
 -\frac{t\dot\varepsilon(t)}{\varepsilon(t)}\longrightarrow1.
 \label{ct:regularization-assumptions}
\end{equation}
Then, for every $x_0,v_0\in\mathcal H$, the solution of \eqref{ct:flow}
satisfies
\begin{equation}
 x(t)\longrightarrow P_Zu\quad\text{strongly},\qquad
 t\|A(x(t))\|\longrightarrow0,\qquad
 t\|\dot x(t)\|\longrightarrow0.
 \label{ct:conclusions}
\end{equation}
Moreover, $t\|V(t,x(t))\|\to0$.
\end{theorem}

\begin{proof}
Throughout the proof, write $a_0=\alpha/2$, $b_0=a_0-1>0$, $p=P_Zu$,
$e=x-p$, and $d=u-p$.  All differential identities below hold almost
everywhere.  Each constant introduced below is independent of $t$ and
retains its assigned value throughout the proof.

\emph{Boundedness.}
Using the first-order representation, introduce
\begin{equation}
 \mathcal H_0(t)=\frac12\|x(t)-p\|^2
 +\frac{a_0}{2b_0}\|w(t)-p\|^2.
 \label{ct:boundedness-energy}
\end{equation}
Differentiation, monotonicity, and Young's inequality give
\begin{align}
 \dot{\mathcal H}_0
 &=-\frac{a_0}{t}\|x-w\|^2-\langle x-p,A(x)\rangle
 -\varepsilon\langle x-p,x-u\rangle
 \notag\\
 &\le-\frac{a_0}{t}\|x-w\|^2-\frac{\varepsilon}{2}\|e\|^2
 +\frac{\varepsilon}{2}\|d\|^2.
 \label{ct:boundedness-estimate}
\end{align}
On the other hand,
\[
 \mathcal H_0
 \le\left(\frac12+\frac{a_0}{b_0}\right)\|e\|^2
 +\frac{a_0}{b_0}\|x-w\|^2.
\]
Choose $T_0\ge t_0$ such that $1/t\ge\varepsilon(t)$ for $t\ge T_0$, and set
$\mu=b_0/(b_0+2a_0)>0$ and $C_{\mathrm H}=\|d\|^2/2$.
The preceding bound gives
$\mu\mathcal H_0\le\frac12\|e\|^2+a_0\|x-w\|^2$, so
\cref{ct:boundedness-estimate} implies
\[
 \dot{\mathcal H}_0
 \le-\varepsilon\bigl(a_0\|x-w\|^2+\tfrac12\|e\|^2\bigr)
       +C_{\mathrm H}\varepsilon
 \le-\mu\varepsilon\mathcal H_0+C_{\mathrm H}\varepsilon.
\]
Multiplication by the integrating factor
$\exp(\mu\int_{T_0}^t\varepsilon(s)\,ds)$ yields, for $t\ge T_0$,
\[
 \mathcal H_0(t)\le\frac{C_{\mathrm H}}{\mu}+
 \left(\mathcal H_0(T_0)-\frac{C_{\mathrm H}}{\mu}\right)
 e^{-\mu\int_{T_0}^t\varepsilon(s)\,ds}
 \le\max\{\mathcal H_0(T_0),C_{\mathrm H}/\mu\}.
\]
Continuity on $[t_0,T_0]$ then gives
$K:=\sup_{t\ge t_0}\mathcal H_0(t)<\infty$.
Thus \cref{ct:boundedness-energy} gives $\|e(t)\|\le\sqrt{2K}$ and bounds
both $x$ and $w$. Set
$M_0=\sqrt{2K}+\|d\|$ so that $\|x(t)-u\|\le M_0$ for all $t\ge t_0$.

\emph{A corrected fast energy.}
We modify the continuous Fast OGDA energy of
\cite[Eq.~(11)]{Bot-Csetnek-Nguyen} by a term that compensates for the
time-dependent regularization; see Remark~\ref{ct:energy-construction}
for its construction. Put $r=\dot x$, $z=A(x)$, and define
\begin{align}
 P&=2a_0e+2tr+tz,\notag\\
 \mathcal E
 &=\frac12\|P\|^2+2a_0b_0\|e\|^2
   +2a_0t\langle e,z\rangle+\frac12t^2\|z\|^2,\notag\\
 \mathcal F
 &=\mathcal E+2a_0t\varepsilon\|e\|^2.
 \label{ct:fast-energy}
\end{align}
Fix the constants $C_-,C_+>0$ specified in
Appendix~\ref{app:continuous-equivalence}. They depend only on $a_0$, and
for every $t\ge T_0$,
\begin{equation}
 C_-\bigl(\|e\|^2+t^2\|r\|^2+t^2\|z\|^2\bigr)
 \le\mathcal F
 \le C_+\bigl(\|e\|^2+t^2\|r\|^2+t^2\|z\|^2\bigr).
 \label{ct:energy-equivalence}
\end{equation}
The proof of this estimate is given in the same appendix.

Define
\begin{equation*}
 \delta(t)=-\frac{t\dot\varepsilon(t)}{\varepsilon(t)}-1,
 \qquad
 \chi(t)=\dot\varepsilon(t)+\frac{a_0}{t}\varepsilon(t)
         =\frac{b_0-\delta(t)}{t}\varepsilon(t).
\end{equation*}
Then $\delta(t)\to0$.  With $g=-\varepsilon r-(x-u)\cdot\chi$,
equation \eqref{ct:expanded-flow} gives
$\dot P=-2b_0r-(2a_0-1)z-t\dot z+2tg$.
Differentiating \eqref{ct:fast-energy} and using $x-u=e-d$ yields
the exact identity
\begin{align}
 \dot{\mathcal F}={}&-4a_0b_0\langle e,z\rangle
 -b_0t\bigl(\|2r+z\|^2+\|z\|^2\bigr)
 -2t^2\langle r,\dot z\rangle\notag\\
 &-4t^2\varepsilon\|r\|^2-2t^2\varepsilon\langle z,r\rangle\notag\\
 &-2a_0(2b_0-\delta)\varepsilon\|e\|^2
   +4a_0(b_0-\delta)\varepsilon\langle e,d\rangle\notag\\
 &-2t^2\chi\langle2r+z,x-u\rangle.
 \label{ct:corrected-energy-identity}
\end{align}
The correction in $\mathcal F$ cancels the mixed term
$-4a_0t\varepsilon\langle e,r\rangle$.

Monotonicity gives $\langle r,\dot z\rangle\ge0$ almost everywhere:
apply monotonicity to $x(t+h)$ and $x(t)$, divide by $h^2$, and pass to the
limit at times where both curves are differentiable.  Also,
\begin{equation*}
 \|2r+z\|^2+\|z\|^2
 \ge\vartheta\bigl(\|r\|^2+\|z\|^2\bigr),
 \qquad \vartheta=3-\sqrt5>0.
\end{equation*}
Choose $T\ge T_0$ sufficiently large that
\begin{equation*}
 |\delta(t)|\le b_0/2,
 \qquad t\varepsilon(t)\le\min\{1,b_0\vartheta/4\}
 \quad(t\ge T).
\end{equation*}
In particular, $b_0/2\le b_0-\delta\le3b_0/2$.  If
$R(t)=\|r(t)\|^2+\|z(t)\|^2$, then
\[
 2t^2\varepsilon|\langle z,r\rangle|
 \le t^2\varepsilon R\le\frac{b_0\vartheta}{4}tR.
\]
Since $|\delta|\le b_0/2$ and
$\chi=(b_0-\delta)\varepsilon/t$, we have
$|\chi|\le3b_0\varepsilon/(2t)$.  Cauchy--Schwarz gives
\[
 \|2r+z\|\le2\|r\|+\|z\|
 \le\sqrt{(2^2+1^2)(\|r\|^2+\|z\|^2)}=\sqrt5\sqrt R.
\]
Using the previously established bound $\|x-u\|\le M_0$, it follows that
\[
 \begin{aligned}
 2t^2|\chi|\,|\langle2r+z,x-u\rangle|
 &\le2t^2|\chi|\,\|2r+z\|\,\|x-u\|\\
 &\le3\sqrt5\,b_0M_0t\varepsilon\sqrt R.
 \end{aligned}
\]
Now apply Young's inequality $XY\le\gamma X^2+Y^2/(4\gamma)$ with
\[
 X=\sqrt{tR},\qquad
 Y=3\sqrt5\,b_0M_0\sqrt t\,\varepsilon,\qquad
 \gamma=\frac{b_0\vartheta}{4}>0.
\]
This yields
\[
 2t^2|\chi|\,|\langle2r+z,x-u\rangle|
 \le\frac{b_0\vartheta}{4}tR
       +\frac{45b_0M_0^2}{\vartheta}t\varepsilon^2.
\]
Set $C_{\mathrm Y}=45b_0M_0^2/\vartheta$. This fixed constant depends
only on $a_0$ and the previously chosen bound $M_0$.
Dropping the remaining nonpositive terms in
\eqref{ct:corrected-energy-identity}, we obtain
\begin{equation}
 \dot{\mathcal F}
 \le-\frac{b_0\vartheta}{2}tR-3a_0b_0\varepsilon\|e\|^2
 +4a_0(b_0-\delta)\varepsilon\langle e,d\rangle
 +C_{\mathrm Y}t\varepsilon^2
 \quad(t\ge T).
 \label{ct:key-estimate}
\end{equation}

\emph{Identification of weak cluster points.}
For $t\ge T$, use $\|e\|\le\sqrt{2K}$,
$0<b_0-\delta\le3b_0/2$, and $t\varepsilon\le1$.
Define the fixed growth constant
\[
 C_{\mathrm G}=6a_0b_0\sqrt{2K}\,\|d\|+C_{\mathrm Y}.
\]
Dropping the negative terms in \eqref{ct:key-estimate} gives
\[
 \dot{\mathcal F}
 \le6a_0b_0\sqrt{2K}\,\|d\|\varepsilon
       +C_{\mathrm Y}t\varepsilon^2
 \le C_{\mathrm G}\varepsilon.
\]
Hence
\begin{equation*}
 \mathcal F(t)
 \le\mathcal F(T)+C_{\mathrm G}\int_T^t\varepsilon(s)\mathrm{d}s
 \le\mathcal F(T)+C_{\mathrm G}\log(t/T).
\end{equation*}
By \eqref{ct:energy-equivalence}, $A(x(t))\to0$.

By the weak--strong closedness stated in Section~\ref{sec:prelim},
every weak cluster point of $x(t)$ belongs to $Z$.
Choose times realizing the limsup of $\langle x(t)-p,u-p\rangle$.
The corresponding bounded sequence admits a weakly convergent subsequence;
its limit lies in $Z$.  Applying \eqref{pre:projection} at that limit gives
\begin{equation}
 \limsup_{t\to\infty}\langle x(t)-p,u-p\rangle\le0.
 \label{ct:projection-limsup}
\end{equation}

\emph{Decay of the corrected energy.}
Let $\omega_0(t)=\max\{\langle e(t),d\rangle,0\}$.
By \eqref{ct:projection-limsup}, $\omega_0(t)\to0$.
Because $0<b_0-\delta\le3b_0/2$, the mixed term in
\eqref{ct:key-estimate} is bounded above by
$6a_0b_0\varepsilon\omega_0$.
Also $tR\ge\varepsilon t^2R$ for $t\varepsilon\le1$.
Set
\[
 \kappa=\frac{\min\{3a_0b_0,b_0\vartheta/2\}}{C_+}>0.
\]
Using the upper bound in \eqref{ct:energy-equivalence}, the negative
terms in \eqref{ct:key-estimate} are bounded above by
$-\kappa\varepsilon\mathcal F$. Thus
\begin{equation}
 \dot{\mathcal F}
 \le-\kappa\varepsilon\mathcal F+\varepsilon\omega(t),
 \qquad
 \omega(t)=6a_0b_0\omega_0(t)+C_{\mathrm Y}t\varepsilon(t)\longrightarrow0,
 \label{ct:scalar-contraction}
\end{equation}
where $\kappa$ depends only on $a_0$.
Since $\int_{t_0}^{\infty}\varepsilon=\infty$,
Lemma~\ref{pre:continuous-contraction}, applied with $F=\mathcal F$,
$\eta=\kappa$, and $b=\omega$, gives $\mathcal F(t)\to0$.
The three conclusions in \eqref{ct:conclusions} follow from
\eqref{ct:energy-equivalence}.  Finally,
\[
 t\|V(t,x(t))\|
 \le t\|A(x(t))\|+t\varepsilon(t)\|x(t)-u\|
 \longrightarrow0.
\]
\end{proof}

\begin{remark}[Construction of the corrected Fast OGDA energy]
\label{ct:energy-construction}
The uncorrected energy $\mathcal E$ in \eqref{ct:fast-energy} is exactly
the Fast OGDA energy $\mathcal E_\lambda$ of
\cite[Eq.~(11)]{Bot-Csetnek-Nguyen} with $\beta\equiv1$,
$\lambda=a_0=\alpha/2$, and reference zero $p$.
We retain the original residual $z=A(x)$ in this energy; monotonicity
then gives $\langle e,z\rangle\ge0$ and
$\langle r,\dot z\rangle\ge0$ almost everywhere.
The additional term $2a_0t\varepsilon\|e\|^2$ is chosen to compensate
for the anchoring terms in \eqref{ct:expanded-flow}.
\end{remark}


The condition $t\varepsilon\to0$ preserves the original energy scale
and makes the remaining error $C_{\mathrm Y}t\varepsilon^2$
equal to $\varepsilon o(1)$. The projection inequality and
$\int^\infty\varepsilon=\infty$ then yield the contraction
\eqref{ct:scalar-contraction}, forcing the entire energy to zero.
This establishes strong selection and vanishing scaled residual and
velocity through \eqref{ct:energy-equivalence}.
The underlying fast quadratic energy is inherited from Fast OGDA;
the present argument couples its compensating correction with
projection selection and control of the time-dependent errors.

\begin{proposition}[Logarithmic regularization and nonexpansive mappings]
\label{ct:logarithmic}
Let $t_0>1$, $c>0$, and $0<q\le1$.  Then
\begin{equation*}
 \varepsilon(t)=\frac{c}{t(\log t)^q}
\end{equation*}
satisfies \eqref{ct:regularization-assumptions}, so the conclusions of
Theorem~\ref{ct:main} hold.  In particular, if
$T\colon\mathcal H\to\mathcal H$ is nonexpansive with nonempty fixed-point
set and $A=I-T$, then
\[
 x(t)\longrightarrow P_{\operatorname{Fix}T}u\quad\text{strongly},
 \qquad \|x(t)-Tx(t)\|=o(t^{-1}).
\]
\end{proposition}

\begin{proof}
Here $t\varepsilon(t)=c/(\log t)^q\to0$,
$-t\dot\varepsilon(t)/\varepsilon(t)=1+q/\log t\to1$,
and $\int_{t_0}^{\infty}\varepsilon=\infty$ for $q\le1$.
For a nonexpansive $T$, the map $I-T$ is monotone and
$2$-Lipschitz, and its zero set is $\operatorname{Fix}T$.
\end{proof}

The proof does not require $\int_{t_0}^{\infty}t\varepsilon(t)^2\,\mathrm{d}t<\infty$.
The error $C_{\mathrm Y}t\varepsilon^2=\varepsilon\,o(1)$ in the proof is absorbed in
\eqref{ct:scalar-contraction}, allowing the full range $0<q\le1$.
Conditions \eqref{ct:regularization-assumptions} are sufficient;
no necessity assertion is made.




%% file: sections/discrete.tex
\section{A second-order accelerated Halpern-type method}\label{sec:discrete}
We now study the {second-order accelerated Halpern-type} method called SAHM associated with
the anchored flow \eqref{ct:flow}.
The convergence proof uses cocoercivity and a separate discrete energy
estimate to control the variation of the historical operator values.
For $k\ge0$ define the historical operator value
\begin{equation*}
 V_k=A(x_k)+\varepsilon_k(x_k-u).
\end{equation*}
Given arbitrary $x_0,x_1\in\HH$, the SAHM iteration is
\begin{equation}\label{dis:algorithm}
 x_{k+1}=x_k+\frac{k}{k+\alpha}(x_k-x_{k-1})
 -\frac{\alpha s}{2(k+\alpha)}V_k
 -\frac{sk}{k+\alpha}(V_k-V_{k-1}),\qquad k\ge1.
\end{equation}
The value $V_{k-1}$ uses $\varepsilon_{k-1}$ and is stored from the
previous step. Recomputing it with $\varepsilon_k$ defines a different
method. After initialization, \eqref{dis:algorithm} requires one new
evaluation of $A$ per update. When $x_0=x_1=u$, the initial evaluation
can be shared. We call the version with logarithmic anchoring \emph{Log-SAHM}. \cref{alg:sahm} is explicit in the equivalent form of \eqref{dis:algorithm}:
Halpern-type describes the fixed anchor and the projection selected
by the method; these coefficients need not form a convex combination.
Setting $\varepsilon_k=0$ recovers Fast KM
\cite{Bot-Nguyen-FastKM}. The operator form \eqref{dis:algorithm}
is therefore also a Tikhonov-regularized Fast KM update.

Next, the following technical lemmas will be used in the main proof. We now indroduce a discrete energy that is a direct analogue of the continuous energy \eqref{ct:fast-energy}.
\subsection{Technical Lemmas}
Write $a=\alpha$, $h_k=x_k-p$, $d_k=x_k-x_{k-1}$, and set
\[
 \begin{gathered}
 b=\frac{3a-2}{2(a-1)},\qquad B=\frac{a(a-2)}2,\qquad
 D=\frac{a(a-2)}{2(a-1)},\\
 F=\frac{(a-2)(3a-2)}{8(a-1)^2},\qquad c_0=\frac{a-2}{a-1}.
 \end{gathered}
\]
We use the Fast KM energy of \cite[Section~3]{Bot-Nguyen-FastKM} at its
parameter $\lambda=a/2$, evaluated at the actual historical values:
\begin{equation*}
 \begin{split}
 E_k={}&\frac12\norm{a h_k+2k d_k+bskV_{k-1}}^2
       +B\norm{h_k}^2\\
 &+Dsk\ip{h_k}{V_{k-1}}+Fs^2k^2\norm{V_{k-1}}^2.
 \end{split}
\end{equation*}
The point $p$ need not satisfy $A(p)+\varepsilon_k(p-u)=0$.
No sign assumption on $\ip{h_k}{V_{k-1}}$ is needed.
Nevertheless, $BF-D^2/4=a(a-2)^2/[8(a-1)]>0$. Pure algebra therefore
gives constants $m,M>0$, independent of $k$, such that
\begin{equation}\label{dis:equiv}
 mS_k\le E_k\le MS_k,\qquad
 S_k=\norm{h_k}^2+k^2\norm{d_k}^2+k^2\norm{V_{k-1}}^2.
\end{equation}
A proof with explicit constants is in \cref{app:equiv}.

\begin{lemma}[Exact energy difference]\label{dis:identitylemma}
Put $h=h_k$, $d=d_{k+1}$, $d_-=d_k$, $v=V_k$, and
$w=V_k-V_{k-1}$. Then
\begin{equation}\label{dis:recurrence}
 (k+a)d=kd_- -\frac{sa}{2}v-skw
\end{equation}
and
\begin{equation}\label{dis:identity}
 \begin{split}
 E_{k+1}-E_k={}&-a(a-2)s\ip h v
 -(a-2)(2k+a+1)\norm d^2\\
 &-\frac{s(a-2)}{a-1}[(3a-2)k+2a^2-1]\ip d v\\
 &-\frac{s^2(a-2)}{2(a-1)}
       [(3a-2)k+a^2+a-1]\norm v^2\\
 &-c_0sk^2\ip{d_-}{w}+\frac{c_0s^2}{2}k^2\norm w^2.
 \end{split}
\end{equation}
\end{lemma}
\begin{proof}
For a short way to organize the expansion, put
$Q_k=kd_k+skV_{k-1}$. The identities
$a^2/2+B=a(a-1)$, $ab+D=2a$, and $b^2/2+F=b$ yield
\[
 E_k=a(a-1)\norm{h_k}^2+2a\ip{h_k}{Q_k}
          +b\norm{Q_k}^2+(2-b)k^2\norm{d_k}^2.
\]
By \eqref{dis:recurrence},
$Q_k=(k+a)d+s(k+a/2)v$ and
$Q_{k+1}=(k+1)(d+sv)$; hence
$Q_{k+1}-Q_k=-(a-1)d-s(a-2)v/2$.
Subtract the four terms using $h_{k+1}=h+d$.
The $\ip h d$ terms cancel. For the last term use
\[
 k^2\norm{d_-}^2
 =\norm{(k+a)d+sa v/2}^2
       +2sk^2\ip{d_-}{w}-s^2k^2\norm w^2.
\]
Collecting $\norm d^2$, $\ip d v$, and $\norm v^2$, and using
$2-b=c_0/2$, gives \eqref{dis:identity}.
No operator inequality is used in this calculation.
\end{proof}

We then record the cocoercivity bound used for each regularized operator.
\begin{lemma}\label{dis:coco}
If $A$ is $\rho$-cocoercive and $\varepsilon>0$, then
$A+\varepsilon\Id$ is $\rho/(1+\rho\varepsilon)$-cocoercive.
\end{lemma}
\begin{proof}
Put $g=Ax-Ay$ and $d=x-y$. Weighted Cauchy--Schwarz and cocoercivity give
\[
 \norm{g+\varepsilon d}^2
 \le (\rho^{-1}+\varepsilon)
       (\rho\norm g^2+\varepsilon\norm d^2)
 \le(\rho^{-1}+\varepsilon)\ip{g+\varepsilon d}{d}.
\]
\end{proof}

We further calculate the difference between successive energies $E_{k+1}$ and $E_k$.
\begin{lemma}[Dissipation]\label{dis:dissipationlemma}
For some integer $K$ and constants $c_1,c_2,c_3,C>0$, the following
estimate holds for every $k\ge K$:
\begin{equation}\label{dis:dissipation}
 \begin{split}
 E_{k+1}-E_k\le{}&
 -d_0\ip{h_k}{V_k}-c_1k\norm{d_{k+1}}^2
 -c_2k\norm{V_k}^2\\
 &-c_3k^2\norm{V_k-V_{k-1}}^2
     +Cr_k\norm{x_{k-1}-u}^2,\qquad d_0=a(a-2)s.
 \end{split}
\end{equation}
\end{lemma}
\begin{proof}
With $z=sv$, the terms in \eqref{dis:identity} involving only $d,v$
are $-(a-2)[k\mathcal Q(d,z)+\mathcal R(d,z)]$, where
\[
 \mathcal Q=2\norm d^2+2b\ip d z+b\norm z^2,\qquad
 \mathcal R=(a+1)\norm d^2+
   \frac{2a^2-1}{a-1}\ip d z+
   \frac{a^2+a-1}{2(a-1)}\norm z^2.
\]
The matrix of $\mathcal Q$ has determinant $b(2-b)>0$ and trace
$2+b$. Thus $\mathcal Q\ge\mu(\norm d^2+\norm z^2)$ with
$\mu=b(2-b)/(2+b)>0$. Write $\gamma=(2a^2-1)/(a-1)$.
Since $\mathcal R\ge-\gamma(\norm d^2+\norm z^2)/2$, for
$k\ge\gamma/\mu$ these terms are bounded above by
$-Lk(\norm d^2+s^2\norm v^2)$, where $L=(a-2)\mu/2$.

Set $\delta_k=\varepsilon_k-\varepsilon_{k-1}$ and
$e_{k-1}=x_{k-1}-u$. Although $w$ compares two different regularized
operators, it can be decomposed as
\[
 w=g_k+\delta_k e_{k-1},\qquad
 g_k=(A+\varepsilon_k\Id)x_k-(A+\varepsilon_k\Id)x_{k-1}.
\]
Choose any $\nu\in(s/2,\rho)$ and put
$\rho_k=\rho/(1+\rho\varepsilon_k)$.
Eventually $\rho_k\ge(\rho+\nu)/2$. By Lemma~\ref{dis:coco} and Young's
inequality,
\begin{align}
 \ip{d_-}{w}
 &\ge \rho_k\norm{w-\delta_k e_{k-1}}^2
                  +\delta_k\ip{d_-}{e_{k-1}}\notag\\
 &\ge \nu\norm w^2-C_\nu\delta_k^2\norm{e_{k-1}}^2
                  +\delta_k\ip{d_-}{e_{k-1}},\label{dis:defect}
\end{align}
where $C_\nu=2\rho\nu/(\rho-\nu)$ suffices. Indeed the Young bound
for the cross term has coefficients
$\rho_k-\nu$ and $\rho_k^2/(\rho_k-\nu)$, leaving the error coefficient
$\rho_k\nu/(\rho_k-\nu)\le C_\nu$.
The last, linear term in $\delta_k$ must still be absorbed.

Applying \eqref{dis:defect} to the remaining terms leaves
$-\beta k^2\norm w^2+Ck^2\delta_k^2\norm{e_{k-1}}^2+T_k$,
where $\beta=c_0s(\nu-s/2)>0$ and
$T_k=-c_0sk^2\delta_k\ip{d_-}{e_{k-1}}$.
Use \eqref{dis:recurrence} again, obtaining
\[
 T_k=-c_0sk(k+a)\delta_k\ip d e
     -\frac{c_0s^2a}{2}k\delta_k\ip v e
     -c_0s^2k^2\delta_k\ip w e,\qquad e=e_{k-1}.
\]
Three applications of Young's inequality yield
\[
 T_k\le\frac L2k\norm d^2+\frac{Ls^2}{2}k\norm v^2
          +\frac\beta2k^2\norm w^2+C k^3\delta_k^2\norm e^2.
\]
For clarity, the respective error coefficients before replacing lower
powers of $k$ by $k^3$ are
$c_0^2s^2k(k+a)^2/(2L)$,
$c_0^2s^2a^2k/(8L)$, and $c_0^2s^4k^2/(2\beta)$.
The assertion follows with $c_1=L/2$, $c_2=Ls^2/2$, and
$c_3=\beta/2$.
\end{proof}
This estimate is established before boundedness of the iterates.
The factor $k^3$ in \eqref{dis:conditions} is the cost of absorbing
the history mismatch into the velocity dissipation.

\subsection{Main result}
\begin{theorem}[Strong selection and fast original residual]\label{dis:main}
Let $A:\HH\to\HH$ be $\rho$-cocoercive with $Z=\zer A\ne\varnothing$.
Let $\alpha>2$, $0<s<2\rho$, and suppose that the positive sequence
$(\varepsilon_k)_{k\ge0}$ satisfies
\begin{equation}\label{dis:conditions}
 \sum_{k=0}^{\infty}\varepsilon_k=\infty,\qquad
 k\varepsilon_k\longrightarrow0,\qquad
 r_k:=k^3(\varepsilon_k-\varepsilon_{k-1})^2=o(\varepsilon_k)\quad(k\ge1).
\end{equation}
For every anchor $u$ and every $x_0,x_1$, the sequence generated by
\eqref{dis:algorithm} satisfies
\begin{equation}\label{dis:conclusions}
 x_k\longrightarrow P_Zu\ \hbox{strongly},\qquad
 k\norm{A(x_k)}\longrightarrow0,\qquad
 k\norm{x_k-x_{k-1}}\longrightarrow0.
\end{equation}
In particular, for nonexpansive $T$, one may take $A=I-T$ and
$0<s<1$.
\end{theorem}
\begin{proof}[Proof of \cref{dis:main}]
Monotonicity and $A(p)=0$ give
\begin{equation}\label{dis:anchor}
 \ip{h_k}{V_k}\ge
 \varepsilon_k\bigl(\norm{h_k}^2-\ip{h_k}{u-p}\bigr).
\end{equation}
Two elementary consequences of \eqref{dis:equiv} and
\eqref{dis:recurrence}, proved in \cref{app:forward}, are
\begin{align}
 E_k&\le M_1\bigl(\norm{h_k}^2+k^2\norm{d_{k+1}}^2+
       k^2\norm{V_k}^2+k^2\norm{V_k-V_{k-1}}^2\bigr),\label{dis:forward}\\
 \norm{x_{k-1}-u}^2&\le C(E_k+1).\label{dis:oldpoint}
\end{align}
First we prove boundedness. Young's inequality in \eqref{dis:anchor}
gives
\[
 -d_0\ip{h_k}{V_k}\le-\frac{d_0}{2}\varepsilon_k\norm{h_k}^2+
 \frac{d_0}{2}\varepsilon_k\norm{u-p}^2.
\]
Take $\eta_1=d_0/(2M_1)$. Since $k\varepsilon_k\to0$, for large $k$
we have $(d_0/2)k\varepsilon_k\le\min(c_1,c_2)$ and
$(d_0/2)\varepsilon_k\le c_3$.
Thus \eqref{dis:forward} permits absorption of
$\eta_1\varepsilon_k E_k$ into the four negative squared terms.
Equations \eqref{dis:dissipation} and \eqref{dis:oldpoint} give
\begin{equation*}
 E_{k+1}\le(1-\eta_1\varepsilon_k+Cr_k)E_k
                      +C_1\varepsilon_k+Cr_k.
\end{equation*}
As $r_k=o(\varepsilon_k)$, eventually
$E_{k+1}\le(1-\eta_1\varepsilon_k/2)E_k+C_2\varepsilon_k$.
The coefficient belongs to $[0,1]$ for large $k$; comparison with
$\max(E_K,2C_2/\eta_1)$ proves boundedness.

By \eqref{dis:equiv}, $(x_k)$ is bounded,
$d_k=O(1/k)$, and $V_{k-1}=O(1/k)$. Consequently
$A(x_k)=V_k-\varepsilon_k(x_k-u)\to0$.
Weak sequential compactness and weak--strong closedness imply that
every weak cluster point of $(x_k)$ belongs to $Z$.
Taking a subsequence realizing the scalar upper limit and then a weakly
convergent subsequence, \eqref{pre:projection} yields
\begin{equation}\label{dis:projectionlim}
 \limsup_{k\to\infty}\ip{h_k}{u-p}\le0.
\end{equation}
This argument uses residual convergence before the projection property.

Return to \eqref{dis:anchor} without the Young relaxation. In the same
way as above, \eqref{dis:forward} absorbs
$\eta_2\varepsilon_kE_k$ for some $\eta_2>0$, giving
\[
 E_{k+1}\le(1-\eta_2\varepsilon_k)E_k+
 d_0\varepsilon_k\ip{h_k}{u-p}+Cr_k(E_k+1).
\]
Set $b_k=d_0\ip{h_k}{u-p}+C(r_k/\varepsilon_k)(E_k+1)$.
Boundedness and \eqref{dis:projectionlim} give $\limsup b_k\le0$.
By Lemma~\ref{pre:contraction}, $E_k\to0$. Its coercivity gives
$x_k\to p$, $k\norm{d_k}\to0$, and $k\norm{V_{k-1}}\to0$.
Shifting the last index also gives $k\norm{V_k}\to0$.
Finally,
\[
 k\norm{A(x_k)}
 \le k\norm{V_k}+k\varepsilon_k\norm{x_k-u}\longrightarrow0.
\]
\end{proof}

\begin{remark}[Inherited energy and relative error absorption]\label{dis:relative}
The quadratic form of $E_k$, including its coefficients, is inherited
from the Fast KM energy at $\lambda=\alpha/2$
\cite[Section~3]{Bot-Nguyen-FastKM}; the residual values are replaced
by the actual historical values $V_{k-1}$. Its algebraic form is
therefore not a separate contribution.
The time-dependent setting requires additional estimates:
$p=P_Zu$ need not be a zero of the mapping
$x\mapsto A(x)+\varepsilon_k(x-u)$, and
$V_k-V_{k-1}$ contains the history mismatch
$(\varepsilon_k-\varepsilon_{k-1})(x_{k-1}-u)$.
The inherited quadratic coercivity in \eqref{dis:equiv} remains valid for these
historical values. Combined with the defect bound \eqref{dis:defect}
and the anchoring inequality \eqref{dis:anchor}, it allows the
inherited energy to prove strong selection.
Specifically, the proof controls the error relative to the anchoring strength:
$r_k/\varepsilon_k\to0$. It does not require $\sum r_k<\infty$.
For \eqref{dis:schedule}, $r_k$ has order
$1/[k(\log k)^{2q}]$ and is nonsummable when $q\le1/2$.
The contraction argument therefore covers schedules excluded by an
absolute summability argument and proves $E_k\to0$, which yields
both projection selection and the original residual estimate.
\end{remark}
\begin{corollary}[Logarithmic schedules]\label{dis:log}
The conclusions of \cref{dis:main} hold for
\begin{equation}\label{dis:schedule}
 \varepsilon_k=\frac{c}{(k+\kappa)[\log(k+\kappa)]^q},
 \qquad c>0,\quad \kappa>1,\quad 0<q\le1.
\end{equation}
\end{corollary}
\begin{proof}
For $f(t)=c/[t(\log t)^q]$,
$f'(t)=-ct^{-2}(\log t)^{-q}(1+q/\log t)$.
The mean value theorem gives
$\varepsilon_k-\varepsilon_{k-1}
=O(k^{-2}(\log k)^{-q})$ and
$r_k/\varepsilon_k=O((\log k)^{-q})\to0$.
Also $k\varepsilon_k\to0$, and the integral test gives
$\sum_k\varepsilon_k=\infty$ exactly when $q\le1$ in this positive
power family.
\end{proof}

%% file: sections/selection.tex
\section{Limit selection analysis}
\label{sel:section}

The projection target solves a secondary optimization problem over
the fixed-point set:
\begin{equation*}
 p=\underset{z\in Z}{\operatorname{argmin}}\,
       \frac12\norm{z-u}^2 \iff 
\ip{u-p}{z-p}\le0\quad \forall z\in Z.
\end{equation*}
Thus convergence to some point of $Z$ and convergence to $p$ answer
different questions, including in finite-dimensional spaces.
To make this distinction precise for ordinary KM, fix a common
relaxation sequence satisfying $0\le\lambda_k\le1$ and
$\sum_k\lambda_k(1-\lambda_k)=\infty$, and define its selection map by
\[
 R_{T,\lambda}(u)=\operatorname*{weak\,lim}_{k\to\infty}
 x_k,\qquad x_0=u.
\]
The standard KM convergence theorem ensures that this map is well
defined \cite[Chapter~5]{Bauschke-Combettes:book}.
Since the relaxation sequence is common to all initial points, the
finite-iterate maps are nonexpansive. Weak lower semicontinuity of the
norm then gives $\norm{R_{T,\lambda}(u)-R_{T,\lambda}(v)}\le\norm{u-v}$. Points in $Z$ generate constant orbits, so $R_{T,\lambda}$ is a
nonexpansive retraction onto $Z$.

Residual decay, convergence to a limit, and selection of a prescribed zero are distinct properties.  In this section, $T:\mathcal H\to\mathcal H$ is
nonexpansive, $A=I-T$, and
\[
 Z=\operatorname{Fix}T=\operatorname{zer}A\ne\varnothing,
 \qquad p=P_Zu,\qquad d=\|u-p\|.
\]
Thus $Z$ is closed and convex.  We use the same reference point $u$ as the
initial point of KM and as the anchor of SAHM. The
normalized quantities below are used only when $d>0$.

Classical KM does have mathematical descriptions of its limit: Fej\'er
monotonicity yields quantitative restrictions, and its limit map is a
nonexpansive retraction.  The issue is that these descriptions need not
identify the limit with $P_Zu$.  The estimates below organize this distinction
using standard projection and retraction geometry
\cite{bauschke_borwein1996,goebel_sedlak2009}; no novelty is claimed for that
classical framework.

\subsection{Selection error and a variational defect}
\label{sel:metrics-subsection}

For a feasible point $z\in Z$, define
\begin{equation}
 \label{sel:metrics}
 B(z)=\frac{\|z-p\|}{d},\qquad
 Q(z)=\frac{\|z-u\|}{d},\qquad
 G(z)=\frac{\|z-u\|^2-d^2}{d^2}=Q(z)^2-1.
\end{equation}
Here $B$ measures selection error, whereas $G$ measures the relative excess
in the secondary objective $\frac12\|z-u\|^2$.  These definitions concern
feasible points; a smaller distance to $u$ at an infeasible iterate is not
evidence of a better solution to this secondary problem.

\begin{proposition}[Projection geometry and normal-cone defect]
\label{sel:geometry}
For every $z\in Z$, $B(z)^2\le G(z)$.  Suppose that, for some $\eta>0$,
\begin{equation}
 \label{sel:margin}
 \langle p-u,z-p\rangle\ge\eta\|z-p\|
 \qquad\text{for every }z\in Z.
\end{equation}
Then
\begin{equation}
 \label{sel:margin-bound}
 G(z)\ge B(z)^2+\frac{2\eta}{d}B(z),\qquad
 B(z)\le\sqrt{(\eta/d)^2+G(z)}-\eta/d
       \le\frac{d}{2\eta}G(z).
\end{equation}
Independently of \eqref{sel:margin}, let
\begin{equation*}
 \nu(z)=\operatorname{dist}(u-z,N_Z(z)),\qquad
 N_Z(z)=\{n:\langle n,y-z\rangle\le0\ \text{for all }y\in Z\}.
\end{equation*}
Then $\|z-p\|\le\nu(z)$, and $\nu(z)=0$ if and only if $z=p$.
\end{proposition}

\begin{proof}
The projection inequality gives $\langle u-p,z-p\rangle\le0$.  Expanding
the square yields the exact identity
\[
 G(z)=B(z)^2+\frac{2\langle p-u,z-p\rangle}{d^2}.
\]
This proves the first assertion and the lower bound in
\eqref{sel:margin-bound}.  Solving the resulting quadratic inequality and
using $\sqrt{a^2+t}-a=t/(\sqrt{a^2+t}+a)\le t/(2a)$ proves the upper bounds.

For the defect estimate, $u-p\in N_Z(p)$.  Given $n\in N_Z(z)$, put
$e=z-p$ and $r=u-z-n$.  Monotonicity of the normal cone gives
$\langle-r-e,e\rangle\ge0$, whence
$\|e\|^2\le-\langle r,e\rangle\le\|r\|\|e\|$.
Taking the infimum over $n$ proves $\|z-p\|\le\nu(z)$. Since the cone $N_Z(z)$ is closed, then  $\nu(z) = 0 \iff u-z\in N_Z(z)$. By the projection characterization
$u-z\in N_Z(z) \iff z=P_Zu = p$, which finally follows that $\nu(z)=0$ if and only if $z=p$.
\end{proof}

Condition \eqref{sel:margin} holds if the closed ball of radius $\eta$
centered at $u-p$ lies in $N_Z(p)$.  Indeed, for $z\ne p$ one may test the
normal-cone inequality with
$u-p+\eta(z-p)/\|z-p\|$.  At a vertex this condition can hold with
$\eta>0$, so excess anchor distance controls selection error linearly.
At a smooth boundary point the infimum of the ratio in
\eqref{sel:margin} can instead be zero.  The condition is geometric and
does not, by itself, imply a rate for an iterative method.

\subsection{Weak limit of KM iteration}
\label{sel:km-subsection}
The standard KM convergence theorem \cite[Chapter~5]{Bauschke-Combettes:book}
gives a weak limit $z\in Z$; in finite dimensions the convergence is in norm.
The following consequences identify
precisely one sense in which this limit remains close to its initialization.

\begin{proposition}[Dissipation and limit proximity]
\label{sel:km-proximity}
Let $z$ be the weak limit of \eqref{intro:km} and $x_0=u$, the sequence of parameters $\{\lambda_k\}$ satisfies $0\le\lambda_k\le1$, $\sum_{k=0}^{\infty}\lambda_k(1-\lambda_k)=\infty$, and set
\[
 \mathcal D=\sum_{k=0}^{\infty}
       \lambda_k(1-\lambda_k)\|Tx_k-x_k\|^2.
\]
Then, for every $q\in Z$,
\begin{equation}
 \label{sel:dissipation}
 \|z-q\|^2+\mathcal D\le\|u-q\|^2.
\end{equation}
In particular,
\begin{equation}
 \label{sel:km-bounds}
 \|z-p\|\le\sqrt{d^2-\mathcal D},\qquad
 d\le\|z-u\|\le d+\sqrt{d^2-\mathcal D}\le2d.
\end{equation}
Moreover, for every $N\ge0$,
\begin{equation}
 \label{sel:km-tail}
 \operatorname{dist}(x_N,Z)\le\|x_N-z\|
       \le2\operatorname{dist}(x_N,Z).
\end{equation}
For a fixed parameter sequence in \eqref{intro:km}, the weak-limit map $R(u)=z$ is a nonexpansive retraction onto $Z$.
\end{proposition}

\begin{proof}
For $q\in Z$, the convex-combination identity and nonexpansiveness imply
\begin{equation}
 \label{sel:fejer}
 \|x_{k+1}-q\|^2\le\|x_k-q\|^2
       -\lambda_k(1-\lambda_k)\|Tx_k-x_k\|^2.
\end{equation}
Summing and using weak lower semicontinuity gives
\eqref{sel:dissipation}.  Taking $q=p$ and applying the triangle inequality
proves \eqref{sel:km-bounds}; its lower bound follows from $z\in Z$. For \eqref{sel:km-tail}, apply the same Fej\'er inequality to the tail and $q=P_Zx_N$.  Weak lower semicontinuity gives
$\|z-q\|\le\|x_N-q\|$, which proves the upper bound by the triangle inequality.  The lower bound again uses $z\in Z$.

Each map $(1-\lambda_k)I+\lambda_kT$ is nonexpansive.  Consequently, two orbits initialized at $u$ and $v$, with the same parameter sequence, satisfy $\|x_k(u)-x_k(v)\|\le\|u-v\|$.  Passing to their weak limits yields $\|R(u)-R(v)\|\le\|u-v\|$.  Points of $Z$ generate constant orbits, so $R|_Z=I$ and $R(\mathcal H)=Z$.
\end{proof}

In particular, ordinary KM satisfies $0\le B(z)\le1$ and $1\le Q(z)\le2$.
These constants are sharp over the class of constant relaxations
$0<\lambda<1$: for $T(t)=|t|$, $u=-1$, and $1/2<\lambda<1$, one step gives
$z=2\lambda-1\in Z=[0,\infty)$, and the orbit remains there.  Thus
$B(z)=2\lambda-1\to1$ and $Q(z)=2\lambda\to2$ as $\lambda\uparrow 1$.
\subsection{Weak limit of the fast KM}
\label{sel:fast-subsection}

For comparison, the unregularized Fast KM iteration is
\begin{equation}
 \label{sel:fast-km}
 x_{k+1}=x_k+\frac{k}{k+\alpha}(x_k-x_{k-1})
 -\frac{s\alpha}{2(k+\alpha)}Ax_k
 -\frac{sk}{k+\alpha}(Ax_k-Ax_{k-1}),\qquad k\ge1.
\end{equation}
We take $\alpha>2$ and $0<s<1$.  Its known convergence theorem gives
$x_k\rightharpoonup z\in Z$ and
$k(x_k-x_{k-1})\to0$, $kAx_k\to0$ in norm
\cite{Bot-Nguyen-FastKM}.  Ordinary KM's Fej\'er bounds and nonexpansive
limit-map argument cannot be transferred to this inertial iteration.

\begin{proposition}
\label{sel:fast-selection}
For \eqref{sel:fast-km}, the residual partial sums converge weakly and
\begin{equation}
 \label{sel:path-selection}
 z=x_1+\frac{x_1-x_0+sAx_0}{\alpha-1}
 -\frac{s(\alpha-2)}{2(\alpha-1)}
       \sum_{k=1}^{\infty}Ax_k.
\end{equation}
The series is understood in this weak, natural-order sense; in finite
dimensions it converges in norm.  Absolute convergence is not asserted.
If $Z=c+L$ is a closed affine subspace, then
\begin{equation}
 \label{sel:affine-fast}
 z=P_Z\!\left(x_0+\frac{\alpha}{\alpha-1}(x_1-x_0)\right).
\end{equation}
For ordinary KM initialized at $u$, the weak limit in the same affine
setting is $P_Zu$.
\end{proposition}

\begin{proof}
Define
$\Xi_k=(\alpha-1)x_k+k(x_k-x_{k-1}+sAx_{k-1})$.
Substitution of \eqref{sel:fast-km} gives the exact identity
\begin{equation}
 \label{sel:path-increment}
 \Xi_{k+1}-\Xi_k=-\frac{s(\alpha-2)}2Ax_k.
\end{equation}
The stated Fast KM convergence properties imply
$\Xi_k\rightharpoonup(\alpha-1)z$.  Telescoping
\eqref{sel:path-increment} proves both convergence of the partial sums
and \eqref{sel:path-selection}.

Suppose now that $Z=c+L$.  For any $h\in L$ and $t\in\mathbb R$,
nonexpansiveness gives
$\|Tx-c-th\|^2\le\|x-c-th\|^2$.
Its validity for arbitrarily large positive and negative $t$ forces
$\langle Ax,h\rangle=0$.  Hence $P_LAx=0$ for every $x$, without requiring
$T$ to be linear.  Projecting \eqref{sel:path-selection} onto $L$ gives
\[
 P_Lz=P_Lx_1+\frac{P_L(x_1-x_0)}{\alpha-1}.
\]
Since $z\in c+L$, this is equivalent to \eqref{sel:affine-fast}.
Projecting the ordinary KM update instead shows
$P_Lx_k=P_Lu$ for every $k$, and therefore its weak limit is $P_Zu$.
\end{proof}

Thus zero initial momentum, $x_1=x_0=u$, gives projection selection for both methods whenever $Z$ is affine.  The same holds for Fast KM if $x_1-x_0$ is a scalar multiple of $Ax_0$.  Arbitrary tangential initial momentum changes the selected point, as \eqref{sel:affine-fast} shows.
\subsection{A nonaffine example with exactly different selected limits}
\label{sel:nonaffine-subsection}

The following example separates selection from residual convergence under
identical initialization.  It also supplies an exact reference for the
numerical comparisons.

\begin{proposition}[A coupled orthant example]
\label{sel:orthant}
In $\mathbb R^2$, let
\begin{equation*}
 Z=\mathbb R_+^2,\qquad D=\{(a,b):a+b\ge0\},\qquad
 T=P_ZP_D,\qquad u=(-1,-4).
\end{equation*}
Then $T$ is nonexpansive, $\operatorname{Fix}T=Z$, and $p=P_Zu=0$.
Ordinary KM with $x_0=u$ and $\lambda=1/2$ has limit $(49/64,0)$.
Fast KM with $x_0=x_1=u$, $\alpha=6$, and $s=1/2$ has limit
\begin{equation}
 \label{sel:orthant-fast-limit}
 z_{\rm F}=\left(\frac{2137007639}{5872025600},0\right)\ne0.
\end{equation}
In contrast, the ordinary Halpern iteration
$h_0=u$, $h_{n+1}=(u+(n+1)Th_n)/(n+2)$ satisfies
\begin{equation}
 \label{sel:orthant-halpern}
 h_n=\frac1{n+1}(2-3\cdot2^{-n},-4),
\end{equation}
thus $
 n\|h_n\|\longrightarrow\sqrt{20},
 n\|Ah_n\|\longrightarrow\sqrt{17}$.
\end{proposition}

\begin{proof}
The composition of the two projections is nonexpansive.  Since $Z\subset D$,
$T(\mathbb R^2)\subset Z$ and $T|_Z=I$, proving the fixed-point assertion.
The assertion about $p$ is immediate.  Three KM updates give
$(49/64,-1/2)$.  Thereafter $T(a,b)=(a,0)$, the first coordinate is constant,
and the second is halved at each update, proving the KM limit.

Write the Fast KM coordinates as $(a_k,b_k)$.  Exact substitution in
\eqref{sel:fast-km} through index $13$ gives
\[
\begin{array}{c|cc}
 k&a_k&b_k\\ \hline
 12&112686810917/259543531520&-13183673/31682560\\
 13&354799859/855638016&-10456591/29245440.
\end{array}
\]
Both rows lie in the branch $a>0$, $b<0$, $a+b>0$, on which $T(a,b)=(a,0)$.
In particular, the current and historical first-coordinate residuals vanish.
On the candidate continuation in this branch,
\[
 \delta_{k+1}=\frac{k}{k+6}\delta_k,\qquad
 \ell=a_k+\frac{k}{5}\delta_k,
 \qquad \delta_k=a_k-a_{k-1}.
\]
The second expression is invariant.  The displayed prefix gives
$\delta_{13}<0$ and $\ell=2137007639/5872025600$, so $a_k$ decreases to
$\ell>0$.  For the negative second coordinate, put $c_k=-b_k$.  Its candidate
continuation satisfies
\[
 c_{k+1}=\frac{3(k+3)c_k-kc_{k-1}}{2(k+6)}.
\]
If $r_k=c_k/c_{k-1}\in[1/2,1]$, then
\[
 \frac{k+9}{2(k+6)}\le r_{k+1} = \frac{3(k+3)-k/r_{k-1}}{2(k+6)}
       \le\frac{2k+9}{2(k+6)}<1.
\]
The initial ratio $r_{13}=135935683/158204076$ lies in this interval.
Thus $b_k$ increases to zero; the upper ratio bound has an infinite product
equal to zero.  Finally,
\[
 \ell+b_{13}=\frac{1911897749}{299473305600}>0.
\]
The entire candidate tail consequently has $a_k\ge\ell>0$, $b_k<0$, and
$a_k+b_k\ge\ell+b_{13}>0$.  The branch is therefore invariant, closing
the induction and proving \eqref{sel:orthant-fast-limit} for the original
nonlinear iteration, rather than merely for a putative linear tail.

For Halpern, the relevant branch is
$b<0$, $a+b<0$, $a-b>0$, where $T(a,b)=((a-b)/2,0)$. The vector in \eqref{sel:orthant-halpern} satisfies these inequalities for
every $n\ge0$, thus
$$ 
Th_n=(\frac{1}{n+1}\frac{(6-3\cdot 2^{-n})}{2},0) = (\frac{{3-3\cdot 2^{-(n+1)}}}{n+1},0),
$$
the next Halpern iterate is
$$ 
\begin{aligned}
      h_{n+1} &= \frac{1}{n+2}(u + (n+1)Th_n) \\
      &= \frac{1}{n+2}(-1,-4) + \frac{n+1}{n+2}(\frac{{3-3\cdot 2^{-(n+1)}}}{n+1},0) \\
      &= \frac{1}{n+2}(2-3\cdot 2^{-(n+1)},-4).
\end{aligned} 
$$
Therefore,~\eqref{sel:orthant-halpern} holds. 
\end{proof}

For the same operator and anchor, Theorem~\ref{dis:main} ensures that the
SAHM iteration converges strongly to $p=0$ with $\|Ax_k\|=o(k^{-1})$.
Thus this example rigorously distinguishes prescribed-point selection from convergence to an unspecified zero, and distinguishes its residual asymptotic
from that of the displayed ordinary Halpern orbit. 

%% file: sections/experiments.tex
\section{Numerical experiments}
\label{num:section}
\subsection{Problems and common initialization}
\label{num:problems}

For the first three problems, let
\begin{equation*}
 D_m=\{x\in\mathbb R^m:\mathbf 1^\top x\geq0\},\qquad
 T=P_ZP_{D_m},\qquad
 P_{D_m}x=x-\min\{\mathbf 1^\top x/m,0\}\mathbf 1.
\end{equation*}
The sets and anchors are given in \cref{num:problem-table}.
In each case $Z\subset D_m$; hence $T(\mathbb R^m)\subset Z$ and $T|_Z=I$,
which proves $\operatorname{Fix}T=Z$.  Projections onto the orthant and the box
are coordinatewise truncations.  Projection onto the triangle is projection
onto the nonnegative simplex ball.  Each map is nonexpansive and nonlinear.

\begin{table}[t]
 \centering
 \caption{The three polyhedral examples.  Every method uses the same initial
 point $u$ within each row; the anchored methods also use $u$ as their anchor.
 In all three examples $P_Zu=0$.}
 \label{num:problem-table}
 \begin{tabular}{llll}
  \toprule
  Name & $m$ & $Z$ & $u$ \\
  \midrule
  Orthant & $2$ & $\mathbb R_+^2$ & $(-1,-4)$ \\
  Box & $20$ & $[0,1]^{20}$ & $u_i=-1-3(i-1)/19$ \\
  Triangle & $2$ & $\{x\geq0:x_1+x_2\leq1\}$ & $(-1,-4)$ \\
  \bottomrule
 \end{tabular}
\end{table}

The fourth problem has a curved fixed-point set.  Set
\begin{equation*}
 \begin{aligned}
 a&=0.3,\qquad b=\sqrt{1-a^2},\qquad u=(0,-1),\\
 B_\pm&=\{x\in\mathbb R^2:\|x-(\pm a,b)\|\leq1\},
 \qquad T=P_{B_+}P_{B_-}.
 \end{aligned}
\end{equation*}
\begin{figure}[!htbp]
  \centering
  \begin{tikzpicture}[scale=2, font=\small, >=Stealth,
    dim/.style={thin, black!55, {Bar}-{Bar}}]
  \def\a{0.3}
  \pgfmathsetmacro{\b}{sqrt(1-\a*\a)}   
  \pgfmathsetmacro{\ang}{asin(\b)}      

  \begin{scope}
    \clip (-\a,\b) circle (1);
    \fill[green!13] (\a,\b) circle (1);         
  \end{scope}
  \draw[thin, dashed, black!55] (-\a,\b) circle (1);
  \draw[thin, dashed, black!55] (\a,\b)  circle (1);

  \draw[thick] (-\a,\b) +({-\ang}:1) arc ({-\ang}:{\ang}:1);
  \draw[thick] (\a,\b)  +({180-\ang}:1) arc ({180-\ang}:{180+\ang}:1);

  \draw[->] (-1.62,0) -- (1.62,0) node[below right] {$x_1$};
  \draw[->] (0,-1.2) -- (0,2.28) node[above] {$x_2$};

  \fill (\a,\b) circle (1.1pt);
  \fill (-\a,\b) circle (1.1pt);
  \node[left,  inner sep=2pt] at (-1.30,\b) {$\partial B_-$};
  \node[right, inner sep=2pt] at ( 1.30,\b) {$\partial B_+$};
  \draw[dim] (0,\b) -- (\a,\b) node[midway, above, inner sep=1.5pt] {$a$};
  \draw[dim] (\a,0) -- (\a,\b) node[midway, right, inner sep=1.5pt] {$b$};

  \draw[dashed, black!60, ->] (0,-1) -- (0,-0.14);
  \fill (0,-1) circle (0.5pt);
  \node[right, inner sep=5.5pt] at (0,-1) {$u=(0,-1)$};
  \node[star, star points=5, star point ratio=2.1, fill=black,
        minimum size=9pt, inner sep=0pt] at (0,0) {};
  \node[below right, inner sep=3pt, yshift=-1pt] at (0,0) {$p=P_Zu=(0,0)$};
\end{tikzpicture}
\caption{The lens $Z=B_-\cap B_+$.}
\label{fig:lens}
\end{figure}
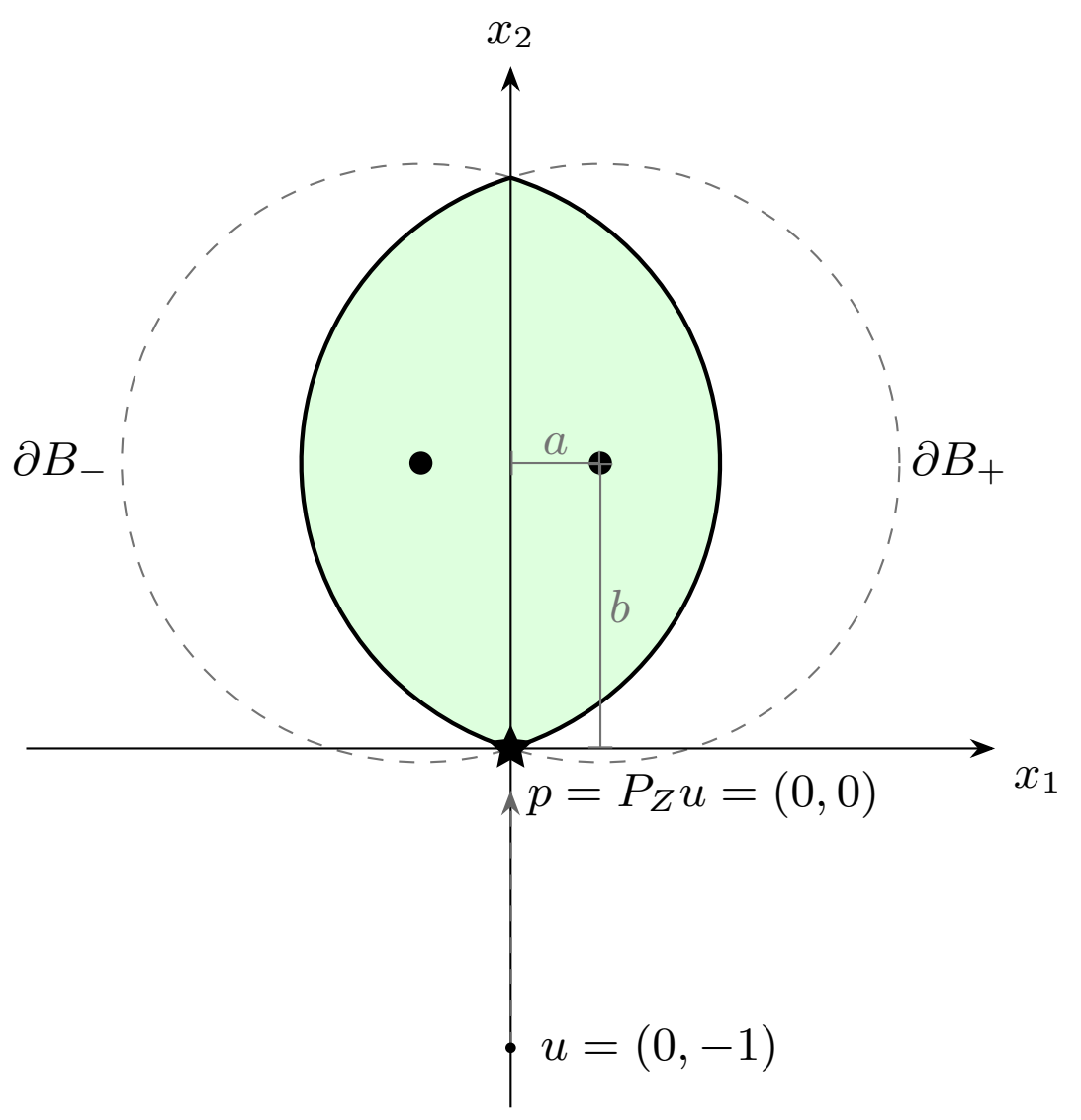
Here $Z=B_-\cap B_+$ is a lens with nonempty interior shown in Figure~\ref{fig:lens}. 

We compare the logarithmic anchoring version of SAHM, denoted by Log-SAHM, with unregularized Fast KM, classical Halpern iteration, the parameter-free adaptive Halpern rule, and classical KM iteration with two relaxation parameters. Log-SAHM combines prescribed-anchor selection with inertial and
operator-difference corrections.  Log-SAHM and the unregularized Fast KM baseline use
\begin{equation}\label{num:parameters}
 \alpha=6,\qquad s=\tfrac12,\qquad x_0=x_1=u,
 \qquad
 \varepsilon_k=\frac{4}{(k+10)\log(k+10)}
\end{equation}
for Log-SAHM, and $\varepsilon_k=0$ for Fast KM. The stored value $V_{k-1}$ uses the historical parameter $\varepsilon_{k-1}$.  The classical baselines are
\begin{equation*}
 \begin{aligned}
 z_0&=u,& z_n&=(1-\lambda)z_{n-1}+\lambda Tz_{n-1},
 &&\lambda\in\{1/2,0.9\},\\
 h_0&=u,& h_n&=\frac{1}{n+1}u+\frac{n}{n+1}Th_{n-1},
 &&n\geq1.
 \end{aligned}
\end{equation*}
The common parameters in \eqref{num:parameters} were chosen during preliminary exploration and then held fixed across the four main examples.  They were not
chosen independently for each row of the performance tables.  All eight configurations, including the two additional regularization coefficients and the adaptive rule used below, have identical recorded initial points within each example.

The adaptive Halpern iteration is the rule
\begin{equation}\label{num:adaptive}
   h_0=u,\quad h_n =\frac{u}{\varphi_n+1}+\frac{\varphi_n}{\varphi_n+1}Th_{n-1},
\end{equation}
where $ \varphi_n=1+\frac{2\langle h_{n-1}-Th_{n-1},\,u-h_{n-1}\rangle}{\|h_{n-1}-Th_{n-1}\|^2}$.

\subsection{Cost and evaluation criteria}
\label{num:protocol}

The budget $n$ counts evaluations of $T$ used to produce an output.  Write
$\widehat x_n$ for this output, regardless of the internal algorithmic index.
For Fast KM and Log-SAHM, $n$ updates produce the internal iterate $x_{n+1}$.
The identical initial values $T(x_0)=T(x_1)$ are reused.  Consequently, all
methods require exactly $n$ solver evaluations to produce $\widehat x_n$.
The final residual requires one additional diagnostic evaluation, treated
equally for all methods. These KM/Halpern comparisons use $N=10^6$
solver evaluations; the additional extragradient comparisons specify
their oracle costs and budgets in \cref{num:flexible,num:flexible-shift}.

With $d=\|u-p\|>0$, we record
\begin{equation*}
 E_n=\frac{\|\widehat x_n-p\|}{d},\qquad
 R_n=\frac{\|\widehat x_n-T\widehat x_n\|}{d},\qquad nR_n.
\end{equation*}
The residual is that of the original operator, not $\|V_k\|$.
For a feasible output, $E_n$ agrees with the selection measure $B$ in
\eqref{sel:metrics}.  When reporting proximity to the anchor through $Q$ and
$G$, we first form $y_n=P_Z\widehat x_n$ and apply those measures to $y_n$;
this avoids interpreting an infeasible point closer to $u$ as an improved
solution of the secondary problem.
To compare the combined requirements, define the finite-horizon statistic
\begin{equation}\label{num:joint-statistic}
 K_{\tau,N}=\min\left\{n\leq N:
       E_j\leq\tau\ \hbox{and}\ R_j\leq\tau
       \text{ for every }j=n,\ldots,N\right\},
\end{equation}
with value $+\infty$ if the set is empty.  This records sustained accuracy only through the observed horizon, not a guarantee about all subsequent iterations.

\subsection{Selected limits versus finite-budget outputs}
\label{num:limits}

\Cref{num:limit-table} reports limits, rather than rounded final iterates.
The orthant selections include the exact limits established in Proposition~\ref{sel:orthant}.
\begin{table}[t]
 \centering\small
 \setlength{\tabcolsep}{4pt}
 \caption{Selected limits from common initial points compared with KM and Fast KM iterations.}
 \label{num:limit-table}
 \begin{tabular}{lccc}
  \toprule
  Set & KM $(\lambda=1/2)$ & KM $(\lambda=0.9)$ & Fast KM \\
  \midrule
  Orthant & $(49/64,0)$ & $(5/4,0)$
    & $(\tfrac{2137007639}{5872025600},0)$ \\
  Box & $q_{B,1/2}$ & $q_{B,0.9}$ & $q_{B,F}$ \\
  Triangle & $(5/8,0)$ & $(4/5,0)$
    & $(\tfrac{27651896071}{87346380800},0)$ \\
  Lens & $(-.06224096,.02185472)$ & $(-.08234703,.02992040)$
    & $(-.02433692,.00799757)$ \\
  \bottomrule
 \end{tabular}
\end{table}

For completeness, the box vectors are listed below.  The displayed nonzero
coordinates are rounded; each zero in the tail is exact, and $0_j$ denotes
$j$ consecutive zero coordinates.  The exact rational values are retained in
the accompanying certificate data.
\begin{equation*}
 \begin{aligned}
 q_{B,1/2}\simeq{}&(.7418914474,.7024177632,.6234703947,\\
                 &\quad .5445230263,.3997861842,.2418914474,.0839967105,0_{13}),\\
 q_{B,0.9}\simeq{}&(.8000000000,.7842105263,.7684210526,\\
                 &\quad .7526315789,.6184210526,.4605263158,.3026315789,.1447368421,0_{12}),\\
 q_{B,F}\simeq{}&(.5444827555,.4645768663,.3653575289,\\
                 &\quad .2375379800,.0822748221,0_{15}).
 \end{aligned}
\end{equation*}

In the lens example, KM with relaxation $1/2$ and $0.9$ enters a single-disk
radial regime after five and two evaluations, respectively.  The subsequent
radial segment remains in the other disk, so its endpoint identifies the
limit.  Its coordinates were evaluated at high precision.  The Fast KM lens
limit was estimated using 80-decimal-digit arithmetic, multiple budgets, and
an algebraically equivalent auxiliary-state recurrence.  The projected
positions at 20000 and 100000 evaluations differed by $2.81\times10^{-17}$;
eight decimal places are reported without claiming a rigorous truncation
error bound.  At $10^6$ evaluations the direct double-precision Fast KM
implementation differs from this estimate by about $9\times10^{-9}$ despite
a residual near $10^{-13}$, illustrating that residuals do not detect
roundoff drift along the fixed-point set.

\begin{figure}[t]
 \centering
 \includegraphics[width=\textwidth]{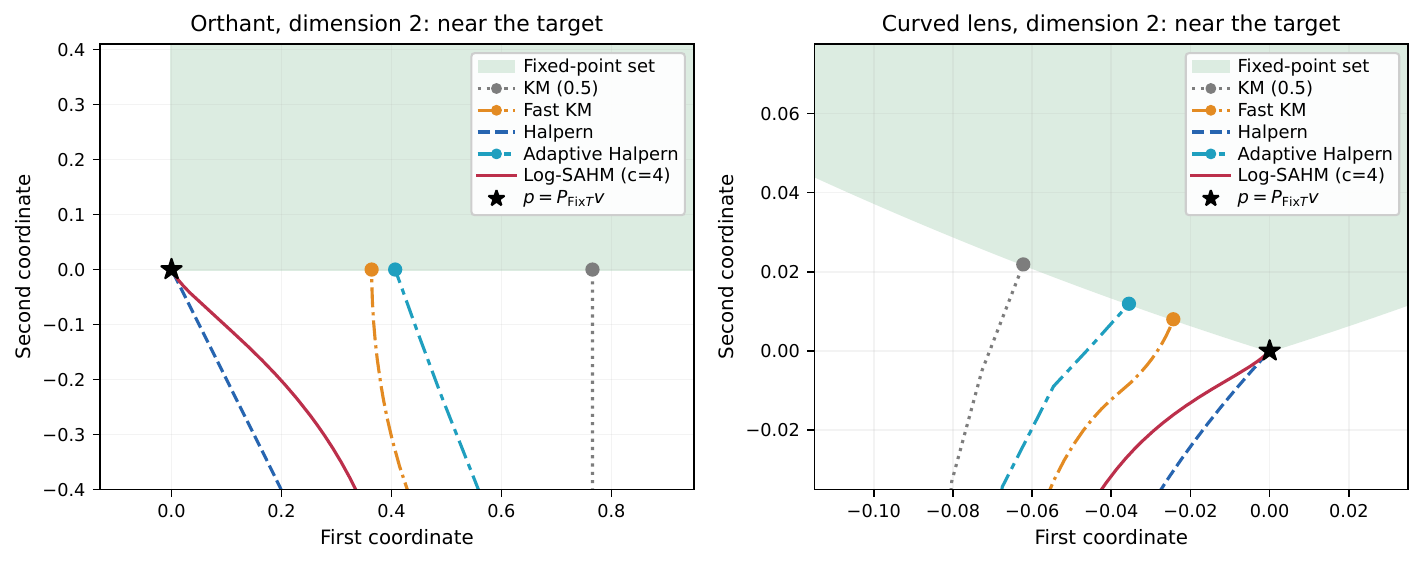}
 \caption{Trajectories near the target in the orthant and lens examples.
 The star marks $p=P_Zu=0$; shaded regions are fixed-point sets.  Line styles
 are identified in the upper-right legend.  Dots mark finite outputs after
 $10^6$ evaluations: for KM and Fast KM they coincide with the limits in
 \cref{num:limit-table}, whereas the adaptive rule stops at a biased fixed
 point.}
 \label{num:geometry-figure}
\end{figure}

The theoretical limit $0$ of an anchored method should likewise not be
confused with its finite output.  For example, at $10^6$ evaluations in the
lens problem, Halpern and Log-SAHM return, respectively,
\[
 (-7.1633\times10^{-7},-7.7472\times10^{-7}),\qquad
 (-2.0740\times10^{-7},-2.2430\times10^{-7}).
\]
These nonzero outputs are consistent with convergence to the same target.
The adaptive rule \eqref{num:adaptive} behaves differently, and its limit
requires no asymptotics beyond the identity below.  In the orthant example the
trajectory enters the branch $T(x_1,x_2)=(x_1,0)$ with $x_2<0$, where
$F(x)=(0,x_2)$ and $\varphi_n+1=-8/x_{2,n}$; the update \eqref{num:adaptive}
then gives
\begin{equation*}
 x_{2,n+1}=\frac{x_{2,n}}{2},\qquad
 \lambda_n=\frac{1}{\varphi_n+1}=\frac{|x_{2,n}|}{8}.
\end{equation*}
The second coordinate halves exactly and
$\sum_n\lambda_n=0.9725<\infty$. The strong-convergence theorem that identifies
the Halpern limit with $P_Zu$ requires $\sum_n\lambda_n=\infty$
\cite{Wittmann}; here that requirement fails, and the limit is accordingly not
the projection.  Evaluated with 60 significant
digits, the iterates converge to $(0.4067914590687645,0)$, at distance
$0.4068$ from $p=0$, while their residual is below $10^{-12}$ after $40$
evaluations.  Thus a machine-precision residual does not certify selection,
even for an anchored method whose weights are inferred from observed
quantities.
\Cref{num:geometry-figure} visualizes the different selections.

\subsection{Accuracy at equal cost and calls to joint accuracy}
\label{num:comparison}

\Cref{num:accuracy-table} compares normalized errors after $10^5$ evaluations.
Log-SAHM reduces both errors relative to classical Halpern in all four
examples.  The target-error ratios lie between $0.3474$ and $0.3481$, and the
residual ratios are approximately $0.3474$.  At $10^6$ evaluations the
corresponding ratios are approximately $0.2895$.  The adaptive rule exhibits the
opposite extreme: it already satisfies $R_n\leq10^{-5}$ at $n=17$ and
$R_n\leq10^{-12}$ at $n=40$ in all four examples, while $E_n$ saturates at
$3.7$--$9.9\times10^{-2}$, so its target errors are about four orders of
magnitude larger than those of Log-SAHM.

\begin{table}[t]
 \centering\small
 \setlength{\tabcolsep}{3.5pt}
 \caption{Normalized errors after $n=10^5$ solver evaluations.  All entries
 are finite-budget measurements.  $E_n$ measures error to the prescribed
 nearest fixed point, rather than error to each method's own limit.  The
 residual of the adaptive rule \eqref{num:adaptive} is machine zero
 ($\leq10^{-15}$) in all four rows at this budget and is therefore omitted.}
 \label{num:accuracy-table}
 \begin{tabular}{lrrrrrr}
  \toprule
  & Fast KM & \multicolumn{2}{c}{Halpern} & Adaptive & \multicolumn{2}{c}{Log-SAHM} \\
  \cmidrule(lr){2-2}\cmidrule(lr){3-4}\cmidrule(lr){5-5}\cmidrule(lr){6-7}
  Set & $E_n$ & $E_n$ & $R_n$ & $E_n$ & $E_n$ & $R_n$ \\
  \midrule
  Orthant & $8.8266\,10^{-2}$ & $1.0846\,10^{-5}$ & $9.9999\,10^{-6}$ & $9.8661\,10^{-2}$ & $3.7686\,10^{-6}$ & $3.4742\,10^{-6}$ \\
  Box     & $7.0764\,10^{-2}$ & $2.6330\,10^{-5}$ & $9.9999\,10^{-6}$ & $8.1088\,10^{-2}$ & $9.1631\,10^{-6}$ & $3.4742\,10^{-6}$ \\
  Triangle& $7.6781\,10^{-2}$ & $1.0846\,10^{-5}$ & $9.9999\,10^{-6}$ & $8.8114\,10^{-2}$ & $3.7686\,10^{-6}$ & $3.4742\,10^{-6}$ \\
  Lens    & $2.5617\,10^{-2}$ & $1.0551\,10^{-5}$ & $9.9999\,10^{-6}$ & $3.7487\,10^{-2}$ & $3.6661\,10^{-6}$ & $3.4742\,10^{-6}$ \\
  \bottomrule
 \end{tabular}
\end{table}

\Cref{num:calls-table} reports \eqref{num:joint-statistic} with
$\tau=10^{-5}$ and $N=10^6$.  Log-SAHM requires between $62.23\%$ and
$64.95\%$ fewer operator evaluations than classical Halpern under this
criterion. The adaptive rule also fails the joint criterion, and for the opposite reason: its residual condition holds after $17$ evaluations, but $E_n$ remains near $10^{-1}$, so $K_{\tau,N}=+\infty$.

\begin{table}[t]
 \centering
 \caption{Calls to joint normalized accuracy $10^{-5}$, sustained through
 $N=10^6$ evaluations.  The percentage is
 $100(1-K_{\tau,N}^{\text{Log-SAHM}}/K_{\tau,N}^{\mathrm{Halpern}})$.}
 \label{num:calls-table}
 \begin{tabular}{lrrr}
  \toprule
  Set & Halpern & Log-SAHM & Reduction \\
  \midrule
  Orthant & $108465$ & $40864$ & $62.33\%$ \\
  Box & $263307$ & $92288$ & $64.95\%$ \\
  Triangle & $108465$ & $40864$ & $62.33\%$ \\
  Lens & $105513$ & $39850$ & $62.23\%$ \\
  \bottomrule
 \end{tabular}
\end{table}

\begin{figure}[p]
 \centering
 \includegraphics[width=\textwidth]{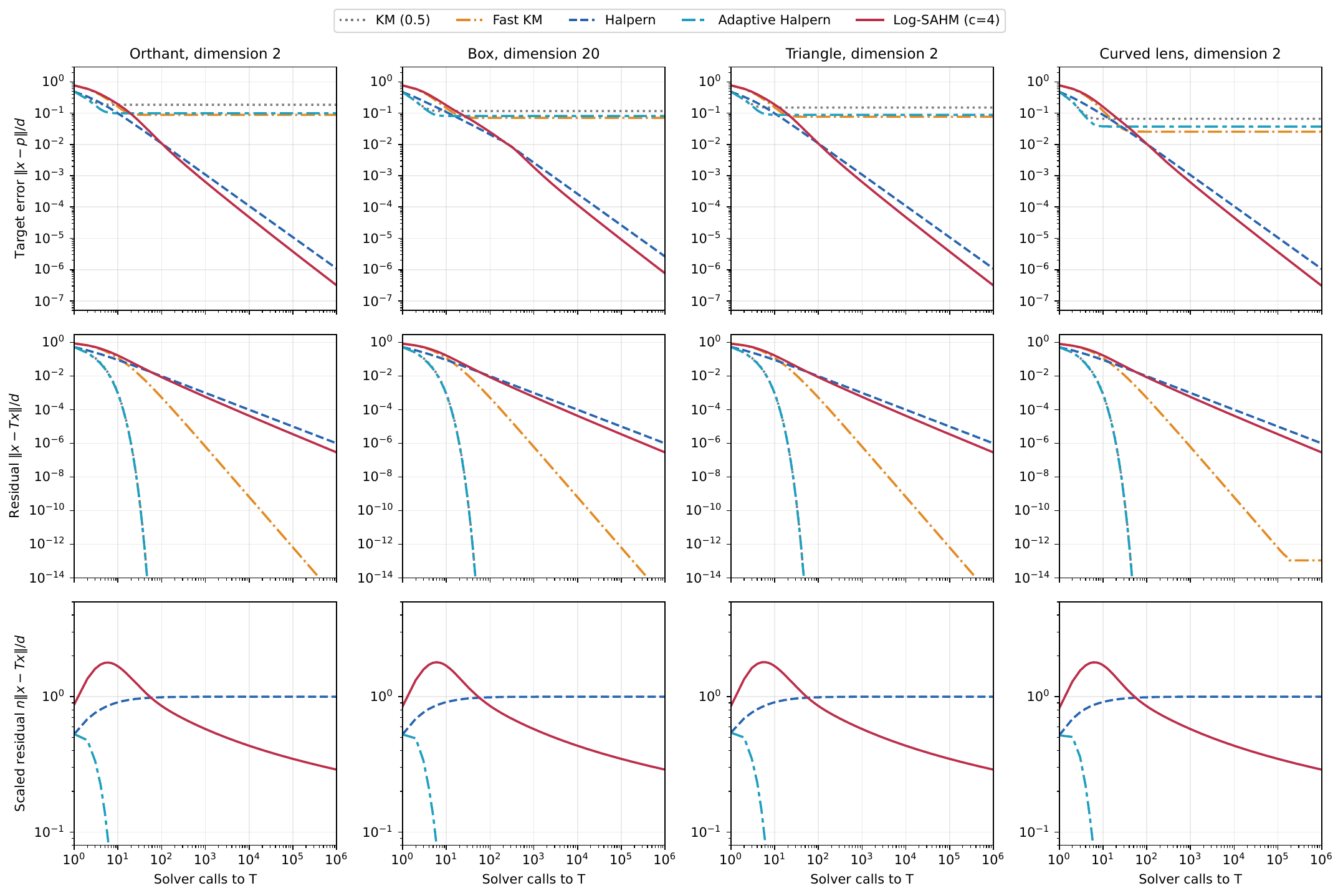}
 \caption{The four main examples through $10^6$ solver evaluations: target
 error $E_n$ (top), original residual $R_n$ (middle), and $nR_n$ for Halpern,
 the adaptive rule, and Log-SAHM (bottom).  KM uses relaxation $1/2$ in the
 plots; its relaxation-$0.9$ limit is included in \cref{num:limit-table}.
 Residuals below $10^{-14}$ are omitted rather than replaced by a positive
 plotting floor. The adaptive rule drives $R_n$ to machine precision within a
 few dozen evaluations while $E_n$ stays on a plateau near $10^{-1}$.}
 \label{num:comparison-figure}
\end{figure}

\subsection{Selection transients, geometry, and limitations}
\label{num:limitations}

The coefficient of regularization affects both selection transients and
late errors.  On the orthant at $10^5$ evaluations, replacing $c=4$ in
\eqref{num:parameters} by $c=1$ gives
$E_n=9.3845\times10^{-3}$, compared with $3.7686\times10^{-6}$ for $c=4$,
although its residual is smaller.  With $c=8$, the target error is
$7.5371\times10^{-6}$.  Thus decreasing the regularization coefficient can
prolong removal of an earlier selection bias, while increasing it can enlarge
the later regularization error.

\begin{figure}[!htbp]
 \centering
 \includegraphics[width=\textwidth]{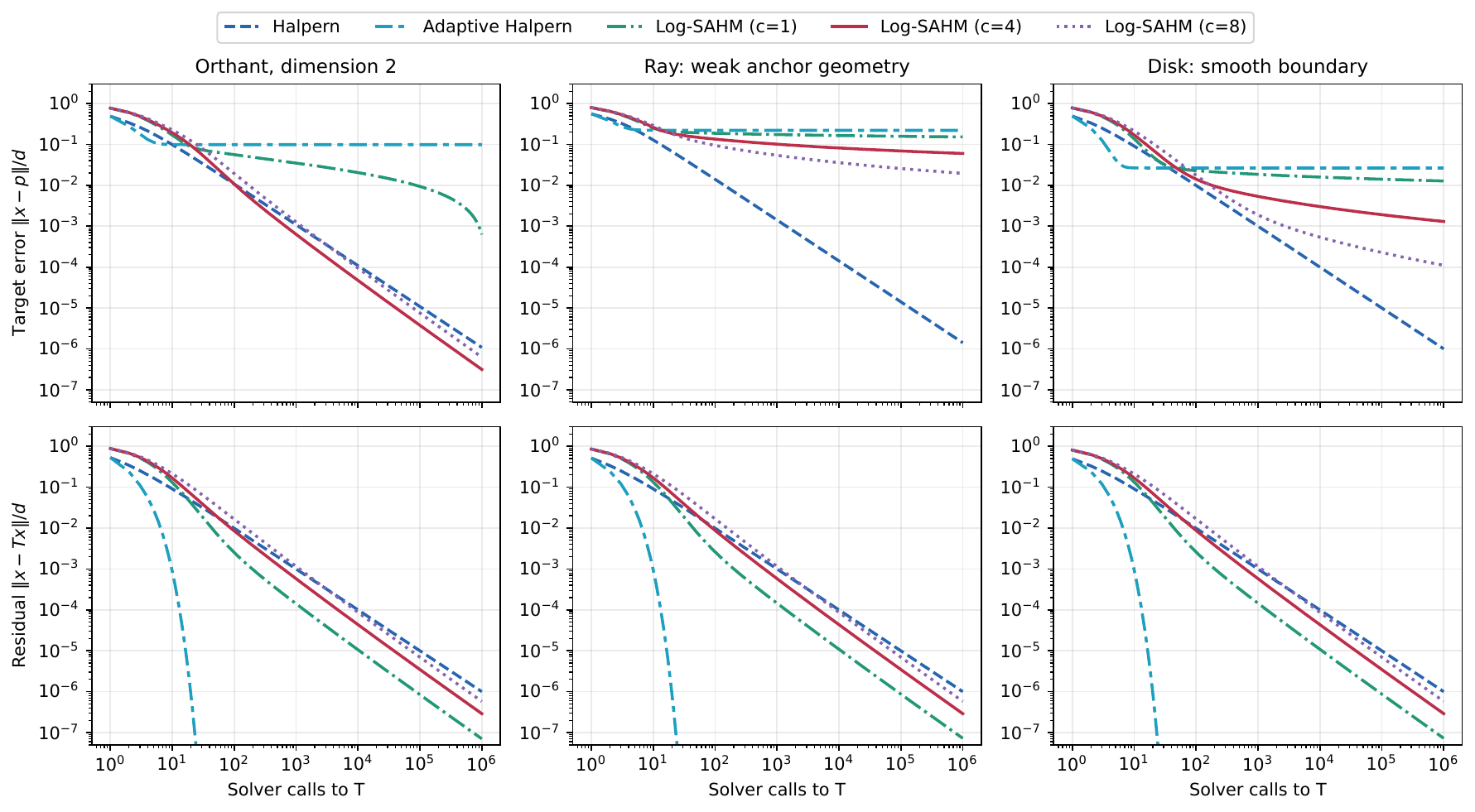}
 \caption{Dependence on the regularization coefficient and on selection
 geometry.  Columns show the orthant, ray, and smooth-disk controls; rows show
 target error and original residual.}
 \label{num:controls-figure}
\end{figure}

The favorable examples satisfy the normal-cone margin condition
\eqref{sel:margin}: one can take $\eta=1$ in the three polyhedral cases and
$\eta=0.3$ in the lens case.  Thus Proposition~\ref{sel:geometry} relates excess anchor
distance to selection error through \eqref{sel:margin-bound}, without itself
implying a distance rate for the iterates.  Nonaffinity by itself is
insufficient.
As a control, let $Z$ be the unit disk, $D=\{x:x_1\leq1\}$,
$T=P_ZP_D$, and $u=(2,2)$.  Then $\operatorname{Fix}T=Z$ and
$p=(1,1)/\sqrt2$.  With the same parameters, at $10^5$ evaluations Halpern
has $E_n\simeq R_n\simeq9.9999\times10^{-6}$, whereas Log-SAHM has
\[
 E_n=1.9079\times10^{-3},\qquad R_n=3.4742\times10^{-6}.
\]
Its target error is approximately $190.8$ times larger despite its smaller
residual.  A further control is
$T(x_1,x_2)=(0,\max\{0,x_2,(x_1+x_2)/2\})$, with $u=(1,0)$ and
$Z=\{0\}\times\mathbb R_+$.  This nonexpansive map also exhibits slow
selection.  Both controls and the coefficient variations appear in
\cref{num:controls-figure}.  They rule out an interpretation of the main
examples as evidence of universal distance acceleration.

The observations support faster approach to the prescribed projection than
classical Halpern on these four constructed examples, while preserving the
distinction between residual decay and selection. 

\input{sections/flexible_comparison}

%% file: sections/flexible_comparison.tex
\subsection{Comparison with anchored extragradient methods}
\label{num:flexible}

We extend the comparison to EG~\cite{Korpelevich1976},
EAG-C and EAG-V~\cite{Yoon-Ryu:21}, FEG~\cite{Lee-Kim}, and
the anchored Popov method (APV)~\cite{Tran-Dinh-Luo}.
The implementations are adapted from the official code. For a prescribed anchor $u$, the exact translation $z=x-u$, with oracle $z\mapsto A(z+u)$, gives the same
recurrences at the required anchor. Every method starts at the same
point; Log-SAHM uses identical initial positions $x_0=x_1$.
All seven previous geometries are retained, with $A=\Id-T$, the shared Lipschitz bound $L=2$, and cocoercivity bound $\rho=1/2$.

The baseline parameters follow the source code: EG uses $0.8/L$,
EAG-C uses $1/(8L)$, EAG-V starts at $0.5/L$, FEG uses $1/L$, and
APV starts at $1/(2\sqrt3 L)$; the variable steps follow the executable
recurrences. In particular, EAG-V uses the denominator
$(k+1)(k+3)$ in the code. Log-SAHM retains
\eqref{num:parameters}, without tuning separately for these comparisons.
We also include the flexible G-EAG scheme of
\cite{Bot-Chenchene-Flexible} with
$\theta=0.8/L$, $\varepsilon_k=1/[\theta(k+5)]$, and a
\emph{G-EAG-log} control with exactly the same logarithmic
$\varepsilon_k$ as Log-SAHM. The latter control separates the effect
of the regularization schedule from the algorithm used with it.

The horizontal axis counts solver evaluations of $A$, each requiring
one evaluation of $T$. For $K\ge1$ completed updates, EG, EAG-C,
EAG-V and both G-EAG variants use $2K$ calls, FEG uses $2K-1$, APV
uses $K+1$, and Log-SAHM uses $K$. We omit FEG's zero-weight first
oracle call and cache APV's previous predictor value. Diagnostic
residual evaluations are excluded consistently. Between completed
updates the last output is held fixed, including when a budget leaves
one unused call. Thus \eqref{num:joint-statistic} remains applicable.

\begin{figure}[!htbp]
 \centering
 \includegraphics[width=\textwidth]{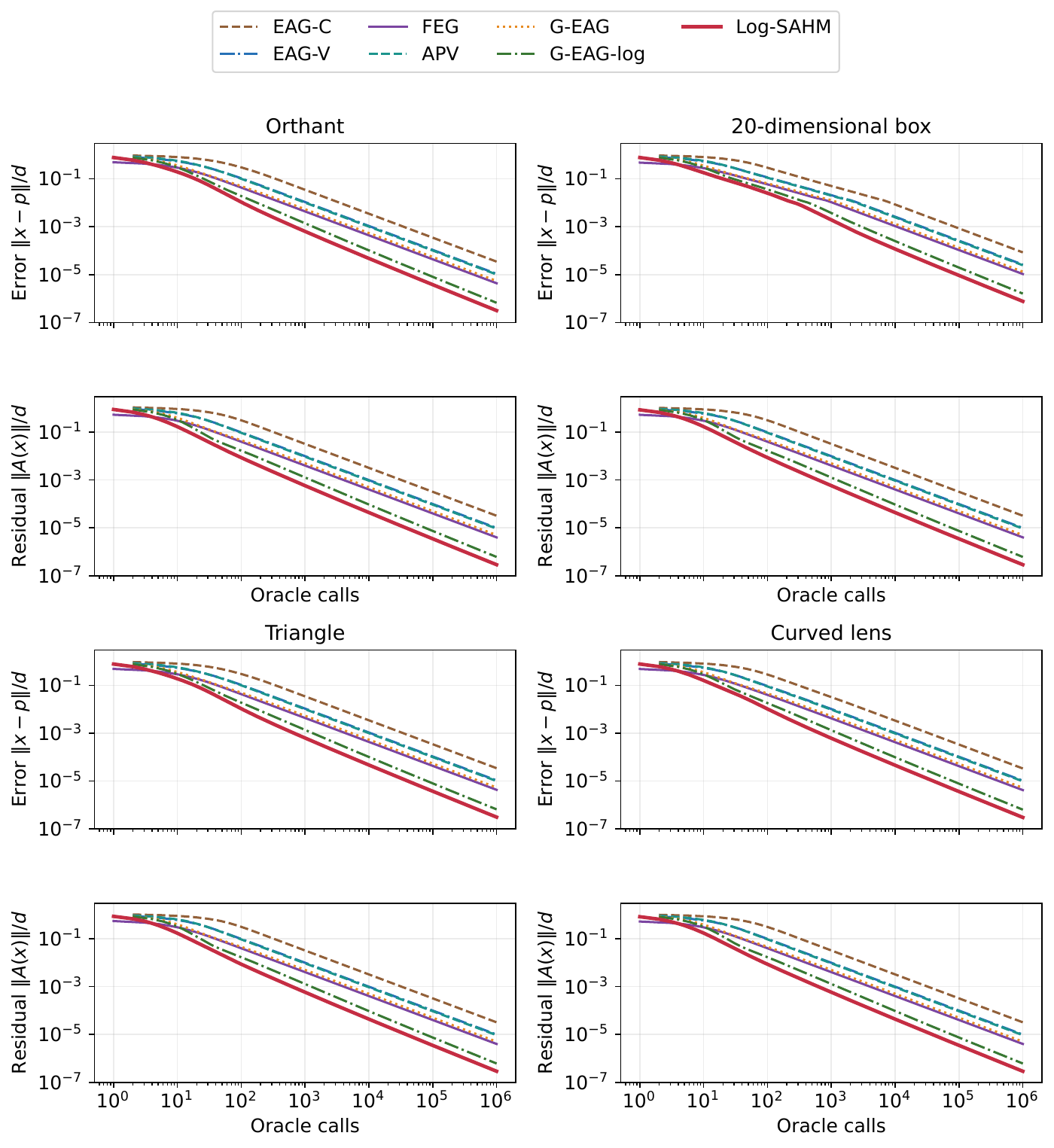}
 \caption{Comparison restricted to anchored methods, from common initial
 points and anchors on the four nonaffine examples.}
 \label{num:flexible-figure}
\end{figure}

\Cref{num:flexible-figure} compares the anchored methods on the four
main geometries; \cref{num:flexible-calls} also retains the unanchored
EG baseline. At $10^6$ calls, Log-SAHM's projection error is approximately
$13.8$ times smaller than FEG's and $2.106$ times smaller than
G-EAG-log's in each case. The joint tolerance $10^{-5}$ is reached
with $90.56\%$--$90.58\%$ fewer calls than FEG on the orthant,
triangle and lens; FEG has not reached that tolerance on the box by
the observation horizon. Against G-EAG-log the reduction is
$49.94\%$--$49.99\%$ across all four cases.
EG reduces the residual faster than the anchored methods in these
examples but retains a nonzero projection error, so its residual
performance alone does not meet the joint task.

The shared-schedule comparison needs a careful interpretation.
After $n$ calls, Log-SAHM has made $n$ updates while G-EAG-log has
made approximately $n/2$. The observed late error ratio is close to
$\varepsilon_{n/2}/\varepsilon_n\simeq2\log n/\log(n/2)$.
These four examples therefore demonstrate oracle efficiency but do
not by themselves isolate a stronger dynamical distance rate.
The one-evaluation SAHM guarantee also uses cocoercivity, whereas
the extragradient baselines address the broader monotone Lipschitz
setting.

\input{sections/flexible_calls_table}

\subsection{The right-shift example}
\label{num:flexible-shift}

To supplement the constructed fixed-point examples, we use the
right-shift problem in \cite{Bot-Chenchene-Flexible} and its code:
\[
 A(x)=x-Sx-(1,-1,0,\ldots),\qquad
 S(x_1,x_2,\ldots)=(0,x_1,x_2,\ldots),\qquad x\in\ell^2.
\]
Its unique zero is $p=(1,0,\ldots)$. Since $S$ is nonexpansive, $A$ is $1/2$-cocoercive and $2$-Lipschitz, so the SAHM theorem applies. We retain the code's common initial point $(0.9,0,1,0,\ldots)$ and
anchor $u=0$. At a budget of $N=2000$ calls, we use $2N+6=4006$ coordinates, sufficient to keep all computed supports and final
residuals away from the truncation boundary. Errors are normalized by $d=\|x_0-p\|=\sqrt{1.01}$.

\begin{figure}[H]
 \centering
 \includegraphics[width=1\textwidth]{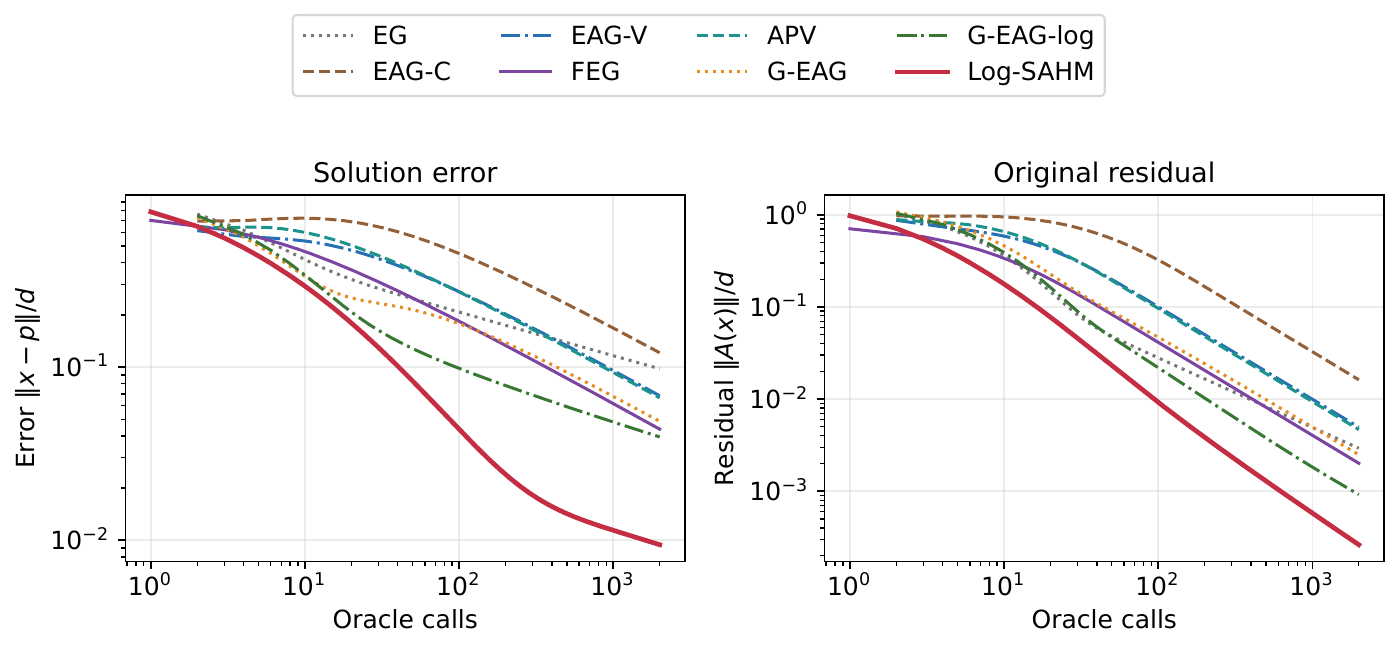}
 \caption{The original right-shift example with common initial points,
 anchors and oracle budgets.}
 \label{num:flexible-controls-figure}
\end{figure}

\Cref{num:flexible-shift-table} compares accuracy under a common
operator-evaluation budget. Its second column records evaluations of $A$.
Under this budget, APV and Log-SAHM complete 1999 and 2000 updates,
respectively, while the other methods complete 1000. These counts
describe how many updates fit within the budget; they are not iteration
counts required to reach a prescribed tolerance. \Cref{num:flexible-shift-table,num:flexible-controls-figure} show
smaller residual and solution error for Log-SAHM in the right-shift
problem. At 2000 calls, its normalized error is $9.4164\times10^{-3}$
and its residual is $2.6263\times10^{-4}$. Relative to FEG these
are smaller by factors $4.65$ and $7.66$, respectively; relative to
G-EAG-log the factors are $4.19$ and $3.50$.
The advantage persists when comparing 1000 updates: Log-SAHM has
error $1.1404\times10^{-2}$ and residual $5.7813\times10^{-4}$,
smaller than FEG's by factors $3.84$ and $3.48$, and smaller than
G-EAG-log's by factors $3.46$ and $1.59$.
This example provides evidence beyond the difference in evaluations
per update, while remaining a finite-horizon observation.

\input{sections/flexible_shift_table}

The geometry controls in \cref{num:limitations}, omitted from
\cref{num:flexible-controls-figure}, delimit the scope of these comparisons.
At $10^6$ calls on
the smooth disk, Log-SAHM has projection error $1.3045\times10^{-3}$,
whereas FEG has $4.0000\times10^{-6}$; on the ray the errors are
$5.9951\times10^{-2}$ and $5.6569\times10^{-6}$, respectively.
Log-SAHM nevertheless has the smaller original residual in both cases. Thus the added baselines support a substantial advantage on
the four principal geometries and the original right-shift problem.

%% file: sections/flexible_calls_table.tex
\begin{table}[t]
 \centering\small
 \caption{Additional baselines: calls to joint normalized accuracy $10^{-5}$ sustained through the $10^6$-call horizon. The output is held fixed between completed updates.}
 \label{num:flexible-calls}
 \begin{tabular}{lrrrr}
  \toprule
  Method & Orthant & Box & Triangle & Lens \\
  \midrule
  EG~\cite{Korpelevich1976} & $>10^6$ & $>10^6$ & $>10^6$ & $>10^6$ \\
  EAG-C~\cite{Yoon-Ryu:21} & $>10^6$ & $>10^6$ & $>10^6$ & $>10^6$ \\
  EAG-V~\cite{Yoon-Ryu:21} & $>10^6$ & $>10^6$ & $>10^6$ & $>10^6$ \\
  FEG~\cite{Lee-Kim} & $433861$ & $>10^6$ & $433861$ & $422053$ \\
  APV~\cite{Tran-Dinh-Luo} & $>10^6$ & $>10^6$ & $>10^6$ & $981005$ \\
  G-EAG~\cite{Bot-Chenchene-Flexible} & $542318$ & $>10^6$ & $542318$ & $527556$ \\
  G-EAG-log & $81710$ & $184338$ & $81710$ & $79676$ \\
  Log-SAHM & $40864$ & $92288$ & $40864$ & $39850$ \\
  \bottomrule
 \end{tabular}
\end{table}

%% file: sections/flexible_shift_table.tex
\begin{table}[t]
 \centering\small
 \caption{Solution error and residual in the original right-shift example under a common budget of 2000 evaluations of $A$. The normalization is $d=\sqrt{1.01}$. FEG's latest completed output costs 1999 evaluations; all other rows use 2000. Smaller errors indicate better accuracy at this budget.}
 \label{num:flexible-shift-table}
 \begin{tabular}{lrrr}
  \toprule
  Method & Evaluations of $A$ & $\|x-p\|/d$ & $\|A(x)\|/d$ \\
  \midrule
  EG~\cite{Korpelevich1976} & $2000$ & $9.7799\times10^{-2}$ & $2.9234\times10^{-3}$ \\
  EAG-C~\cite{Yoon-Ryu:21} & $2000$ & $1.2098\times10^{-1}$ & $1.6180\times10^{-2}$ \\
  EAG-V~\cite{Yoon-Ryu:21} & $2000$ & $6.8288\times10^{-2}$ & $4.9559\times10^{-3}$ \\
  FEG~\cite{Lee-Kim} & $1999$ & $4.3807\times10^{-2}$ & $2.0105\times10^{-3}$ \\
  APV~\cite{Tran-Dinh-Luo} & $2000$ & $6.6414\times10^{-2}$ & $4.6672\times10^{-3}$ \\
  G-EAG~\cite{Bot-Chenchene-Flexible} & $2000$ & $4.8569\times10^{-2}$ & $2.4905\times10^{-3}$ \\
  G-EAG-log & $2000$ & $3.9494\times10^{-2}$ & $9.1991\times10^{-4}$ \\
  Log-SAHM & $2000$ & $9.4164\times10^{-3}$ & $2.6263\times10^{-4}$ \\
  \bottomrule
 \end{tabular}
\end{table}

%% file: sections/discussion.tex
\section{Conlusion}\label{sec:discussion}
The SAHM method combines Halpern-type anchoring with the inertial and
residual-difference structure of Fast KM, giving strong convergence to
$P_Zu$ and $o(k^{-1})$ original residual decay. The logarithmic schedule
supplies infinite cumulative anchoring while its contribution to the
scaled residual vanishes. The history mismatch satisfies
$k^3(\varepsilon_k-\varepsilon_{k-1})^2=o(\varepsilon_k)$ and is absorbed
without summability. The geometry of the fixed-point set determines how
residual decay and selection translate into distance to the prescribed solution.

There is also a cost. Our theorem gives no uniform rate for $\norm{x_k-P_Zu}$, no finite-budget guarantee of superiority over ordinary KM, and no guarantee for the endpoint stepsize $s=2\rho$. When the nearest point lies at a geometrically sharp part of the boundary, residual and selection errors can decrease together,
as in the principal experiments. On a smooth boundary, a tangential selection error can remain after feasibility has essentially been achieved. The control experiments show that this distinction can reverse a finite-budget ranking. In practical problems $P_Zu$ is usually unknown, so the projection error reported here is a benchmark diagnostic rather than a generally implementable stopping criterion.

Natural further questions are quantitative selection rates under error bounds, parameter rules that exploit local geometry, and preconditioned versions preserving the desired projection metric. The established results provide a basis for these questions by making the original residual, the selected solution, and the cost of oracle evaluations explicit.

%% file: sections/energy_appendix.tex
\section{Uniform energy bounds and the forward estimates}\label{app:energy}
This appendix supplies the algebraic details behind
\eqref{ct:energy-equivalence}, \eqref{dis:equiv}, \eqref{dis:forward},
and \eqref{dis:oldpoint}.

\subsection{Continuous energy equivalence}\label{app:continuous-equivalence}
Use the notation of \cref{ct:section}, with $a_0=\alpha/2>1$ and
$b_0=a_0-1$, and put $y=tr$, $\zeta=tz$, and
$S=\|e\|^2+\|y\|^2+\|\zeta\|^2$.
The definition \eqref{ct:fast-energy} becomes
\[
 \begin{aligned}
 P&=2a_0e+2y+\zeta,\\
 \mathcal F&=\tfrac12\|P\|^2+2a_0b_0\|e\|^2
 +2a_0\langle e,\zeta\rangle+\tfrac12\|\zeta\|^2
 +2a_0t\varepsilon\|e\|^2.
 \end{aligned}
\]
Since $A(p)=0$, monotonicity gives
$\langle e,\zeta\rangle=t\langle x-p,A(x)-A(p)\rangle\ge0$.
Consequently,
\[
 \mathcal F\ge\tfrac12\|P\|^2+2a_0b_0\|e\|^2
                         +\tfrac12\|\zeta\|^2.
\]
The identity $2y=P-2a_0e-\zeta$ and Cauchy--Schwarz yield
$\|y\|^2\le\tfrac34\|P\|^2+3a_0^2\|e\|^2
+\tfrac34\|\zeta\|^2$. Hence
\[
 S\le\tfrac34\|P\|^2+(1+3a_0^2)\|e\|^2
       +\tfrac74\|\zeta\|^2\le J\mathcal F,
 \qquad
 J=\max\left\{\tfrac72,\frac{1+3a_0^2}{2a_0b_0}\right\}.
\]
For the upper bound, let $t\ge T_0$, where $T_0$ is the threshold fixed
in the proof of \cref{ct:main}; then $t\varepsilon(t)\le1$.
Using $\|P\|^2\le12a_0^2\|e\|^2+12\|y\|^2+3\|\zeta\|^2$
and $2\langle e,\zeta\rangle\le\|e\|^2+\|\zeta\|^2$ gives
\[
 \mathcal F\le(8a_0^2+a_0)\|e\|^2+6\|y\|^2
                   +(a_0+2)\|\zeta\|^2
 \le(8a_0^2+a_0)S,
\]
where the last inequality uses $a_0>1$.
Thus \eqref{ct:energy-equivalence} holds with $C_-=J^{-1}$ and
$C_+=8a_0^2+a_0$, both depending only on $a_0$.

\subsection{Discrete coercivity and equivalence}\label{app:equiv}
Use the constants in \cref{sec:discrete} and put
$h=h_k$, $y=kd_k$, $z=skV_{k-1}$, and
$\zeta=ah+2y+bz$. Then
\[
 E_k=\tfrac12\norm\zeta^2+P(h,z),\qquad
 P(h,z)=B\norm h^2+D\ip h z+F\norm z^2.
\]
Completing the square gives
\begin{equation*}
 P(h,z)=B\left\|h+\frac{z}{2(a-1)}\right\|^2
                       +\frac{a-2}{4(a-1)}\norm z^2.
\end{equation*}
More quantitatively, let
\[
 \Delta=BF-D^2/4=\frac{a(a-2)^2}{8(a-1)},\quad
 \tau=B+F,\quad \ell=\Delta/\tau.
\]
The positive definite matrix
$\left(\begin{smallmatrix}B&D/2\\D/2&F\end{smallmatrix}\right)$
has smallest eigenvalue at least $\ell$ and largest eigenvalue at most
$\tau$. These two-dimensional scalar matrix bounds apply to
Hilbert-space pairs by the same orthogonal change of coordinates.
Hence
\[
 \ell(\norm h^2+\norm z^2)\le P(h,z)
                       \le\tau(\norm h^2+\norm z^2).
\]
Set $N=\norm h^2+\norm y^2+\norm z^2$. Since $0<b<2<a$ and
$2y=\zeta-ah-bz$,
\[
 N\le\tfrac34\norm\zeta^2+
        (1+3a^2/4)(\norm h^2+\norm z^2)\le J E_k,\quad
 J=\max\left\{\tfrac32,\frac{1+3a^2/4}{\ell}\right\}.
\]
Conversely, $\norm\zeta^2\le3a^2N$, so that
$E_k\le U N$, where $U=3a^2/2+\tau$.
For $S_k$ in \eqref{dis:equiv},
$\min(1,s^2)S_k\le N\le\max(1,s^2)S_k$.
Thus one possible choice is
\[
 m=\frac{\min(1,s^2)}{J},\qquad M=U\max(1,s^2).
\]
All constants are independent of $k$, the regularization schedule,
and the iterates. No monotonicity or iteration identity was used.

\subsection{Forward control and the old anchor displacement}\label{app:forward}
For $k\ge1$, \eqref{dis:recurrence} implies
\[
 k^2\norm{d_k}^2\le3(a+1)^2k^2\norm{d_{k+1}}^2
       +\frac{3s^2a^2}{4}k^2\norm{V_k}^2
       +3s^2k^2\norm{V_k-V_{k-1}}^2.
\]
Also
$k^2\norm{V_{k-1}}^2\le
2k^2\norm{V_k}^2+2k^2\norm{V_k-V_{k-1}}^2$.
Substitution in the upper bound of \eqref{dis:equiv} proves
\eqref{dis:forward}, for example with
\[
 M_1=M\max\{1,\,3(a+1)^2,\,3s^2a^2/4+2,\,3s^2+2\}.
\]
Finally, $x_{k-1}-u=h_k-d_k+(p-u)$ gives
\[
 \norm{x_{k-1}-u}^2
 \le3(\norm{h_k}^2+\norm{d_k}^2+\norm{p-u}^2)
 \le\frac3m E_k+3\norm{p-u}^2.
\]
This proves \eqref{dis:oldpoint}. 